\documentclass[11pt,reqno]{amsart}
\usepackage{t1enc}
\usepackage[latin1]{inputenc}
\usepackage[english]{babel}
\usepackage{amsmath,amsthm}
\usepackage{amsfonts}
\usepackage{latexsym}
\usepackage[dvips]{graphicx}
\usepackage{graphicx}
\usepackage{float}
\usepackage[natural]{xcolor}
\usepackage{algorithm}
\usepackage{algorithmic}
\usepackage{enumerate}
\usepackage{appendix}
\usepackage{multirow}
\usepackage{xcolor} 
\usepackage[colorlinks,linkcolor=blue]{hyperref}
\usepackage{lineno}
\usepackage{setspace}
\usepackage{comment}
\usepackage{fullpage}
\usepackage{amssymb} 
\usepackage{listings}

\usepackage[inline,shortlabels]{enumitem}

\newtheorem{theorem}{Theorem}\numberwithin{theorem}{section}
\newtheorem{lemma}[theorem]{Lemma}\numberwithin{theorem}{section}
\newtheorem{fact}[theorem]{Fact}
\newtheorem{problem}[theorem]{Problem}

\newtheorem{claim}[theorem]{Claim}

\theoremstyle{definition}
\newtheorem{Definition}[theorem]{Definition}
\theoremstyle{plain}

\def\int{\textrm{int}}
\numberwithin{equation}{section}

\begin{document}

\onehalfspacing
\title{
Tur\'{a}n numbers for suspensions of  bipartite 2-regular graphs}

\author{Haoxing Wu}
\address{HW. School of Mathematics and Statistics, Beijing Institute of Technology, Beijing, 100081, China\\
Email: \texttt{(HW) hxingwu2003@163.com}}
\author{Junxue Zhang}
\address{JZ. School of Science, Hebei University of Technology, Tianjin, 300401, China, Email: \texttt{(JZ) jxuezhang@163.com}.}

\date{}

\begin{abstract}
For a graph $R$, let $\widehat R:=K_1+R$ denote its suspension, obtained by adding a new vertex adjacent to every vertex of $R$.
For every fixed nonempty
bipartite $2$-regular graph $R$, we determine
$\operatorname{ex}(n,\widehat{R})$ for all sufficiently
large $n$ and characterize all extremal graphs. Three distinct cases arise, according as $R$ is $2C_6$, a disjoint union of copies of $C_4$, or neither of these. This extends the known exact  Tur\'an results on odd wheels to suspensions of arbitrary bipartite 2-regular graphs.


\medskip
\noindent\textbf{Keywords:} Tur\'{a}n number, suspension, bipartite  graph, extremal graph

\end{abstract}

\maketitle

\section{Introduction}

Throughout this paper, all graphs are finite and simple. 
Let $H$ be a graph. We say that a graph is \emph{$H$-free} if it contains no copy of $H$ as a subgraph. For a fixed graph $H$, the \emph{Tur\'an number} of $H$, denoted by $\operatorname{ex}(n,H)$, is the maximum number of edges in an $n$-vertex $H$-free graph. An $n$-vertex graph $G$ is called an \emph{extremal graph} for $H$ if $G$ is $H$-free and $e(G)=\operatorname{ex}(n,H)$.  
One of the fundamental results in extremal graph theory is the Erd\H{o}s-Stone-Simonovits theorem~\cite{ErdosSimonovits1966}, which states that for every fixed graph $H$ with $\chi(H)\ge 3$,  $\operatorname{ex}(n,H) = \left(1-\frac{1}{\chi(H)-1}\right)\binom{n}{2}+o(n^2). $ 
Thus the chromatic number  determines the asymptotic behavior of the Tur\'an number of  every  non-bipartite graph. However, determining the exact value of $\operatorname{ex}(n,H)$ and characterizing the corresponding extremal graphs are usually much more difficult.



Let $T(n,r)$ be the balanced complete $r$-partite graph on $n$
vertices, and let $t(n,r)=e(T(n,r))$. 
For graphs $G$ and $H$ on disjoint vertex sets, let $G+H$ denote their \emph{join}, obtained by adding all edges between $V(G)$ and $V(H)$ to their disjoint union. 
The \emph{suspension} of a graph $R$ is $\widehat R:=K_1+R$, where the vertex in $K_1$ is usually called the \emph{apex}. 
Tur\'an problems for suspensions have been studied since the early development of extremal graph theory.
Simonovits \cite{simonovits1968method}  determined $\operatorname{ex}\bigl(n,\widehat{T(rd,r)}\bigr)$ for all fixed $r$ and $d$ and all sufficiently large $n$. A graph $H$ is called \emph{edge-critical} if it contains an edge $e$ such that $\chi(H-e)<\chi(H)$. Another fundamental result  of Simonovits \cite{Simonovits1974} states that if $H$ is edge-critical with
$\chi(H)=r+1\ge 3$, then, for all sufficiently large $n$, 
$\operatorname{ex}(n,H)=t(n,r)$,
and $T(n,r)$ is the unique extremal graph.
 In particular, if $H$ is edge-critical, then so is its suspension
$\widehat H$.  Hence the  Tur\'an number and the extremal graphs
of $\widehat H$ are determined for all sufficiently large $n$.  
Beyond the suspension of a single edge-critical graph,
Chen, Gould, Pfender and
Wei~\cite{ChenGouldPfenderWei2003}  determined the Tur\'an numbers and extremal graphs for $\widehat{kK_r}$. 
More recently, Hou, Li and Zeng~\cite{HouLiZeng2024} generalized
this result to 
suspensions of collections of  edge-critical graphs with the same chromatic number. 

\subsection{Suspensions of sparse bipartite graphs}
For sparse bipartite base graphs beyond this general framework, several exact results are known. Erd\H{o}s, F\"uredi, Gould and
Gunderson~\cite{ErdosFurediGouldGunderson1995} determined the
extremal graphs for $\widehat{kK_2}$.
Zhu, Wang, Zhang and Zhang~\cite{ZhuWangZhangZhang2026} 
studied Tur\'an numbers of suspensions of balanced trees, where a tree is called balanced if the sizes of its two color classes differ by at most one.
Among sparse bipartite base graphs, even cycles form a particularly natural family. The suspension of an even cycle is precisely an odd wheel:
$
\widehat{C_{2k}}=W_{2k+1}.
$ For $k=2$, Simonovits' result~\cite{simonovits1968method}
already determines $\operatorname{ex}(n,\widehat{C_4})$
for all sufficiently large $n$.
Dzido and Jastrz\k{e}bski~\cite{DzidoJastrzebski2018} determined $\operatorname{ex}(n,\widehat{C_4})$ and $\operatorname{ex}(n,\widehat{C_6})$. Yuan~\cite{Yuan2021}
subsequently determined the Tur\'an numbers and  all extremal graphs for $\widehat{C_{2k}}$ for every
fixed $k\ge2$ and all sufficiently large $n$.
We summarize these results below.


\begin{theorem}[Simonovits~\cite{simonovits1968method}, Dzido--Jastrz\k{e}bski~\cite{DzidoJastrzebski2018},
Yuan~\cite{Yuan2021}]
\label{thm:yuan-odd-wheel}
Let $k\ge2$. For all sufficiently large $n$,
\[
\operatorname{ex}(n, \widehat{C_{2k}})
=
\begin{cases}
\displaystyle
\left\lfloor\frac{n^2}{4}\right\rfloor
+\left\lfloor\frac n4\right\rfloor
+\left\lfloor\frac{n+1}{4}\right\rfloor,
&\text{if }k=2,\\[1em]
\displaystyle
\max_{n_0+n_1=n}
\left\{
n_0n_1+
\left\lfloor\frac{(k-1)n_0}{2}\right\rfloor+1
\right\},
&\text{if }k\ge3.
\end{cases}
\]
\end{theorem}

\subsection{Our results}

We extend Theorem~\ref{thm:yuan-odd-wheel} from a single even cycle
to an arbitrary nonempty bipartite $2$-regular graph,
that is, a nonempty disjoint union of even cycles. Throughout the remainder of the paper, write
$R\cong \bigcup_{i=1}^{r} C_{2p_i},
2\le p_1\le\cdots\le p_r,
s:=\sum_{i=1}^{r}p_i$.


\begin{theorem}[Main result]\label{thm:main}
Let $R$ be a nonempty bipartite $2$-regular graph on $2s$ vertices. 
Define
\[
f(n,R):=\max\{g(n_0,n_1,R):
n_0,n_1\in\mathbb Z^+,\,n_0+n_1=n\}.
\]
Then, for all sufficiently large $n$,
$\operatorname{ex}(n,\widehat R)=f(n,R)$. 
Moreover, the extremal graphs are precisely those in $\mathcal K_n(R)$.
Here $g(n_0,n_1,R)$ and $\mathcal K_n(R)$ are defined below.
\end{theorem}

The function $g(n_0,n_1,R)$ appearing in the definition of
$f(n,R)$ is given as follows:
\[
g(n_0,n_1,R):=
\begin{cases}
n_0n_1
+6\left\lfloor\frac{n_0}{4}\right\rfloor
+\binom{n_0\bmod 4}{2}
+6\left\lfloor\frac{n_1}{4}\right\rfloor
+\binom{n_1\bmod 4}{2},
& \text{if } R\cong2C_6,\\[1.5em]
\displaystyle
n_0n_1
+\left\lfloor\frac{(s-1)n_0}{2}\right\rfloor
+\left\lfloor\frac{n_1}{2}\right\rfloor,
& \text{if $s$ is even and } R\cong \frac{s}{2}C_4,\\[1.5em]
\displaystyle
n_0n_1
+\left\lfloor\frac{(s-1)n_0}{2}\right\rfloor
+\min\left\{1,\left\lfloor\frac{n_1}{2}\right\rfloor\right\},
& \text{otherwise.}
\end{cases}
\]

We now define the extremal family $\mathcal K_n(R)$.
The extremal graphs have a common underlying structure:
a nearly balanced complete bipartite graph with additional
edges inside its parts. 
For an integer $d\ge1$, a graph $G$ is called \emph{nearly $d$-regular} if every vertex has degree $d$, except possibly one vertex of degree $d-1$. 
For $R=C_{2k}$ with $k\ge3$, the known extremal graphs
are obtained by inserting a nearly $(k-1)$-regular
$P_{2k-1}$-free graph into one part and a single edge
into the other. Our constructions retain this form in the general case,
but differ when $R\cong rC_4$ or $R\cong 2C_6$.
Let $R^-:=\{R-x:x\in V(R)\}$. 
 Every member of $R^-$ is a path together with a (possibly empty) collection of cycles.
For a graph family $\mathcal F$, a graph $G$ is called $\mathcal F$-free if $G$ contains no $F\in\mathcal F$ as a subgraph. We regard $K_0$ as the empty graph.

\begin{Definition}[Extremal graphs]
For $n_0,n_1\in\mathbb{Z}^+$, let $\mathcal K_{n_0,n_1}(R)$ be the
family of all graphs $G_0+G_1$, where $|V(G_i)|=n_i$ for
$i\in\{0,1\}$, and
$$
\begin{cases}
\begin{aligned}
G_0\cong
\left\lfloor\frac{n_0}{4}\right\rfloor K_4
\cup K_{n_0\bmod4},
\qquad
&G_1\cong
\left\lfloor\frac{n_1}{4}\right\rfloor K_4
\cup K_{n_1\bmod4},
\end{aligned}
&\hspace{-0.5em} \text{if }R\cong2C_6,\\[1.5em]
\begin{aligned}
G_0\text{ is nearly $(s-1)$-regular and $R^-$-free},
\quad
&G_1\text{ is nearly $1$-regular},
\end{aligned}
&\hspace{-0.5em} \text{if $s$ is even and } R\cong \frac{s}{2}C_4,\\[1.5em]
\begin{aligned}
G_0\text{ is nearly $(s-1)$-regular and $R^-$-free},
\quad
&e(G_1)=\min\left\{1,\left\lfloor\frac{n_1}{2}\right\rfloor\right\},
\end{aligned}
&\hspace{-0.5em} \text{otherwise.}
\end{cases}
$$

Let $\mathcal K_n(R)$ denote the graphs with the maximum number of edges among $\bigcup_{\substack{n_0+n_1=n\\n_0,n_1\in\mathbb{Z}^+}}
\mathcal K_{n_0,n_1}(R)$.
\end{Definition}

Our result concerns suspensions of disjoint unions of even cycles,
where every vertex of the cycles is adjacent to a common apex.
Related Tur\'an problems for disjoint unions of odd wheels,
each with its own apex, have also been studied. Xiao and Zamora~\cite{xiao2021note}
determined the Tur\'an numbers and extremal graphs for
disjoint unions of identical odd wheels, and conjectured
a corresponding result for odd wheels of different orders.
A graph $F$ is called \emph{vertex-critical} if there exists
a vertex $v\in V(F)$ such that
$\chi(F-v)<\chi(F)$.
Zhang~\cite{zhang2026extremal} subsequently extended their result
to disjoint unions of vertex-critical graphs with the
same chromatic number under a suitable ordering condition,
thereby resolving the conjecture of Xiao and Zamora. 
 In particular, Zhang's result implies that, for any fixed vertex-critical graph $H$ with $\chi(H) \geq 3$, the Tur\'an problem for $tH$ is determined by that for $H$. Since every suspension $\widehat R$ considered here is vertex-critical, Theorem~\ref{thm:main}, together with Zhang's result, determines $\operatorname{ex}(n,t\widehat R)$ for every fixed positive integer $t$ and all sufficiently large $n$.


\subsection{Proof outline and organization}
We briefly outline the proof of Theorem~\ref{thm:main}.
The constructions in $\mathcal K_n(R)$ give the lower bound $f(n,R)$. 
For the upper bound and the characterization of extremal graphs,
we use Simonovits' progressive induction
method~\cite{simonovits1968method}, following its application
to odd wheels by Yuan~\cite{Yuan2021}.
We apply this method to
$\phi(n):=\operatorname{ex}(n,\widehat R)-f(n,R)$.
Starting from an extremal $\widehat R$-free graph,
we find a large complete bipartite subgraph $K_{N,N}$.
This yields either $\phi(n)<\phi(n-2N)$,
or $\phi(n)<\phi(n-1)$, or an almost balanced partition
of bounded internal degree.
In the last case, local extremal lemmas determine the edge bound
and the equality cases, completing the proof.

The remainder of the paper is organized as follows.
Section~\ref{sec:notation} states lemmas ensuring the existence
of the extremal constructions and the key local extremal lemmas.
Section~\ref{sec:proof-main} uses these lemmas to prove
Theorem~\ref{thm:main} by progressive induction.
Section~\ref{sec:proof-lemmas} proves the local extremal lemmas.  Finally, in Section~\ref{sec:concluding}, we pose an open
problem concerning the extremal graphs for $tK_1+R$.

\section{Preliminaries and key lemmas}\label{sec:notation}

\subsection{Basic notation}

For a graph $G$ and a vertex $x\in V(G)$, let $N_G(x)$ be the neighborhood of $x$ in $G$, and let $d_G(x):=|N_G(x)|$. 
If $X,Y\subseteq V(G)$ are disjoint, then $G[X,Y]$ denotes the bipartite subgraph of $G$ with vertex classes $X$ and $Y$ and edge set $\{xy\in E(G):x\in X,\ y\in Y\}$. Let $e_G(X,Y):=e(G[X,Y])$.
For $x\in V(G)$ and $X\subseteq V(G)$, let $N_{G[X]}(x):=N_G(x)\cap X$ and $d_{G[X]}(x):=|N_{G[X]}(x)|$. If $x\in X$, then let $N_{G[X]}[x]:= \left\{x\right\} \cup N_{G[X]}(x)$. For $\emptyset\ne S\subseteq V(G)$, set $N_{G[X]}(S):=\bigcap_{v\in S}N_{G[X]}(v)$, and for a nonempty
subgraph $G'\subseteq G$, write
$N_{G[X]}(G'):=N_{G[X]}(V(G'))$.
For a graph $H$ and an integer $q\ge0$, we write $qH$ (or $q\cdot H$) for the disjoint union of $q$ copies of $H$. 
For $m\ge1$, let $P_m$ denote the path on $m$ vertices and let $S_m:=K_{1,m-1}$. 
We use the convention that
$\binom{a}{b}=0$ when $b>a$.

Throughout the paper, all $O(1)$ terms are bounded independently of $n,n_0,n_1$,
with implicit constants depending only on the fixed graph $R$.
For sufficiently large $n$, every maximizing pair $(a,b)$ satisfies
$|a-b|=O(1)$.
Indeed, for fixed $R$,
$
g(n_0,n_1,R)=n_0n_1+\alpha n_0+\beta n_1+O(1)
$
for some constants $\alpha,\beta$. Write
$a=n/2+x$ and $b=n/2-x$. Since  
$g(a,b,R)\ge g(\lfloor n/2\rfloor,\lceil n/2\rceil,R)$, we obtain
$x^2\le |\alpha-\beta||x|+O(1)$.
Hence $|x|=O(1)$, and therefore $|a-b|=O(1)$.
We may choose a maximizing pair satisfying $a\ge b$.

\subsection{Properties of the construction family}\label{subsec:construction}

\begin{lemma}\label{lem:U-nonempty}
Let $R\cong\bigcup_{i=1}^r C_{2p_i}$ with $r\ge 1$ and $2\le p_1\le\cdots\le p_r$, and let $s=\sum_{i=1}^r p_i$.
Then, for all sufficiently large $n$, there exists an $n$-vertex nearly $(s-1)$-regular $R^-$-free graph. 
Moreover, when $r\ge2$, this graph may be chosen to be $\left\{C_{2p_1},\ldots,C_{2p_r}\right\}$-free.
\end{lemma}

When $R\cong2C_6$, the nonemptiness of $\mathcal K_{n_0,n_1}(R)$ is immediate. In the other two cases,
it follows from Lemma~\ref{lem:U-nonempty}. 
Thus $\mathcal K_{n_0,n_1}(R)\neq\emptyset$ whenever $n_0$ is sufficiently large, and every graph in $\mathcal K_{n_0,n_1}(R)$ then has $g(n_0,n_1,R)$ edges. Consequently, for all sufficiently large $n$, $\mathcal K_n(R)\neq\emptyset$,
and every graph in $\mathcal K_n(R)$ has exactly $f(n,R)$ edges.
The following lemma shows that the graphs in $\mathcal K_{n_0,n_1}(R)$ are indeed $\widehat R$-free.

\begin{lemma}\label{lem:construction-free}
For all $n_0,n_1\in\mathbb{Z}^+$, every graph
$G\in\mathcal K_{n_0,n_1}(R)$ is
$\widehat{R}$-free. Consequently,
$\operatorname{ex}(n,\widehat{R})\ge \max\left\{g(n_0,n_1,R):n_0+n_1=n\right\}= f(n,R)$ for all sufficiently large $n$.
\end{lemma}

\begin{lemma}\label{lem:key_construction}
For any $N_0$,
there exists a constant $N>N_0$, divisible by $4$,
such that, for all sufficiently large $n$, there exists a graph $H_n\in\mathcal K_n(R)$ with a partition $V(H_n)=X_{n_0}\cup Y_{n_1}$ satisfying $H_n[X_{n_0},Y_{n_1}]\cong K_{n_0,n_1}$ and
$\max\{\Delta(H_n[X_{n_0}]),\Delta(H_n[Y_{n_1}])\}\le s-1$.
Moreover, there exist sets $B'_1\subseteq X_{n_0}$ and
$B'_2\subseteq Y_{n_1}$, each of size $N$, such that
\begin{center}
\begin{enumerate*}[label=\textup{(\roman*)}, itemjoin=\qquad]
\item\label{1} $H_n[B'_1,B'_2]\cong K_{N,N}$;
\item\label{2} $e\big(H_n-(B'_1\cup B'_2)\big)=f(n-2N,R)$;
\item\label{3} $e(H_n[B'_1,X_{n_0} \setminus B'_1])=e(H_n[B'_2,Y_{n_1} \setminus B'_2])=0$;
\item\label{4} 
$g(N,N,R)\le
e(H_n[B'_1\cup B'_2])+1$.
\end{enumerate*}
\end{center}
\end{lemma}

\subsection{Key lemmas for the main theorem}


\begin{theorem}[Erd\H{o}s--Stone~\cite{ErdosStone1946}]\label{thm:erdos-stone}
For all integers $p\ge2$ and $N\ge1$, and every $\eta>0$, there exists an integer $n_0$ such that every graph on $n\ge n_0$ vertices with at least $t(n,p-1)+\eta n^2$ edges contains $T(pN,p)$ as a subgraph.
\end{theorem}

\begin{theorem}[Erd\H{o}s--Sachs~\cite{ErdosSachs1963}]\label{thm:erdos-sachs}
For integers $k\ge 2$ and $\ell\ge 3$, there exists a simple $k$-regular graph on at most $4\sum_{i=1}^{\ell-2}(k-1)^i$
vertices with girth at least $\ell$.
\end{theorem}

\begin{lemma}\label{lem:f-f}
Let $R$ be a bipartite $2$-regular graph. Then there exists a constant
$h=h(R)$, depending only on $R$, such that, for all sufficiently large $n$, $f(n,R)-f(n-1,R)\ge n/2+h(R)$.
\end{lemma}

\begin{lemma}\label{lem:special-two-C6}
Let $G$ be an $n$-vertex $\widehat{2C_6}$-free graph with a partition
$V(G)=V_0\cup V_1$. For $i\in\{0,1\}$, write
$G_i:=G[V_i]$, $n_i:=|V_i|$, and $k_i:=\Delta(G_i)\le5$.
Suppose $n_0\ge n_1$. Then, for $n_0$ sufficiently large,
$e(G)\le f(n,2C_6)$, with equality if and only if
$G\in\mathcal K_n(2C_6)$.
\end{lemma}

\begin{lemma}\label{lem:special-all-C4}
Let $r \in \mathbb{Z}^+$ and let $G$ be an $n$-vertex $\widehat{rC_4}$-free graph with a partition $V(G)=V_0\cup V_1$. For $i\in\{0,1\}$, write
$G_i:=G[V_i]$, $n_i:=|V_i|$, and $k_i:=\Delta(G_i)\le2r-1$.
Suppose $n_0\ge n_1\ge 2r^2+r$. Then, for $n_0$ sufficiently large,
$e(G)\le g(n_0,n_1,rC_4)$, with equality if and only if $G\in\mathcal{K}_{n_0,n_1}({rC_4})$.
\end{lemma}

\begin{lemma}\label{lem:bipartite}
Let $R\cong \bigcup_{i=1}^r C_{2p_i}$ with $2\le p_1\le\cdots\le p_r$ and $s=\sum_{i=1}^r p_i$. Suppose $R \not\cong rC_4$ and $R\not\cong2C_6$. Let $G$ be an $n$-vertex $\widehat{R}$-free graph with a partition $V(G)=V_0\cup V_1$.
For $i\in\{0,1\}$, write
$G_i:=G[V_i]$, $n_i:=|V_i|$, and $k_i:=\Delta(G_i)\le s-1$.
Suppose $n_0\ge n_1\ge 8s$.
Then, for $n_0$ sufficiently large, 
$e(G)\le g(n_0,n_1,R)$, with equality if and only if $G\in\mathcal{K}_{n_0,n_1}(R)$.
\end{lemma}

\section{Proof of Theorem \ref{thm:main}}\label{sec:proof-main}

\begin{proof} 
Fix $R$ and choose $N_*=N_*(R)$ sufficiently large so that  Lemmas~\ref{lem:special-two-C6},~\ref{lem:special-all-C4}, and~\ref{lem:bipartite} apply whenever $n_0=n_1\ge N_*$. Then apply Lemma~\ref{lem:key_construction} with $N_0=N_*$ and let
$N>N_*$ be the resulting constant, which is divisible by $4$.
By Lemma~\ref{lem:U-nonempty} and the fact that every maximizing pair
satisfies $|n_0-n_1|=O(1)$, there exists an integer $M_0$ such that
$\mathcal K_n(R)\neq\emptyset$ for every $n\ge M_0$.
Moreover, Lemma~\ref{lem:construction-free} yields $\operatorname{ex}(n,\widehat R)\ge f(n,R)$ whenever $n\ge M_0$.
For each $n\ge M_0$, let $\mathcal A_n$ be the set of extremal $\widehat R$-free graphs and define $\phi(n):=\operatorname{ex}(n,\widehat R)-f(n,R)\ge0$.
For an extremal graph $\Gamma$, we write
$\phi(\Gamma):=\phi(|V(\Gamma)|)$.
The proof uses a progressive induction argument formalized in  Claim~\ref{Claim:mainthm}. 

\begin{claim}\label{Claim:mainthm}
Suppose that there exists $M_1>M_0+2N$ such that for every $n>M_1$ and $L_n\in\mathcal A_n$, at least one of the following holds: 
\begin{center}
\begin{enumerate*}[label=\textup{(\alph*)}, itemjoin=\qquad]
\item\label{a} $\phi(n-2N)\ge\phi(n)+1$;
\item\label{b} $\phi(n-1)\ge\phi(n)+1$;
\item\label{c} $L_n\in\mathcal K_n(R)$.
\end{enumerate*}
\end{center}
Then $L_n\in\mathcal K_n(R)$ for every $L_n\in\mathcal A_n$ and all sufficiently large $n$.
\end{claim}
\begin{proof}


Let 
$\Phi_0:=\max\{\phi(m):M_1-2N<m\le M_1\}$.
Let $L_n\in\mathcal A_n$.
Suppose $L_n\notin\mathcal K_n(R)$ for some $n>M_1+2N\Phi_0$. 
Set $\Gamma_0:=L_n$. 
Since $\Gamma_0$ does not satisfy \ref{c},
either \ref{a} or \ref{b} applies, yielding an extremal graph
$\Gamma_1$ with $\phi(\Gamma_1)>\phi(\Gamma_0)\ge0$.
Since any graph satisfying \ref{c} has $\phi=0$, none of the subsequent graphs
can satisfy \ref{c}. Repeating this process, we obtain $\Gamma_1,\ldots,\Gamma_t$ such that
$1\le |V(\Gamma_{i-1})|-|V(\Gamma_i)|
\le 2N$ and
$\phi(\Gamma_i)\ge\phi(\Gamma_{i-1})+1$ for $i\in[t]$. We stop at $\Gamma_t$ once $|V(\Gamma_t)|\le M_1$.
So
$M_1-2N<|V(\Gamma_t)|\le M_1$.
Thus
$\Phi_0\ge\phi(\Gamma_t)\ge
\phi(\Gamma_0)+t \ge
t$ and $2N\Phi_0<{n-M_1}\le{n-|V(\Gamma_t)|}\le 2Nt$, 
a contradiction. Then
$L_n\in\mathcal K_n(R)$ for sufficiently large $n$.
\end{proof}

By Claim \ref{Claim:mainthm}, it remains to prove that at least one of \ref{a}, \ref{b}, and \ref{c} holds.
Since
$e(L_n)\ge f(n,R)>
\max_{n_0+n_1=n}\{n_0n_1+1\}\ge n^2/4$, 
Theorem~\ref{thm:erdos-stone} with $p=2$ and
$\eta=1/4$ yields, for all sufficiently large $n$,  a copy of
$T(2N,2)=K_{N,N}$ in $L_n$. Let $B_1$ and $B_2$ denote its two parts, each of size $N$. 
For each $i\in [2]$, we claim that $\Delta(L_n[B_i])\le s-1$. If not, then some vertex of $B_i$ would have at least $s$ neighbors in $B_i$. Together with any $s$ vertices of $B_{3-i}$, we obtain a copy of $K_1+K_{s,s}$, which contains $\widehat R$, a contradiction.
Set $A:=V(L_n)\setminus(B_1\cup B_2)$. 
Let $\varepsilon:=2s/N$. For $i\in[2]$, define $C_i:=\{x \in A: d_{L_n[B_i]}(x)<s \text{ and } d_{L_n[B_{3-i}]}(x)\ge(1-\varepsilon)N \}$, and let 
$D:=A\setminus(C_1\cup C_2)$. Since $N$ is sufficiently large,
$(1-\varepsilon)N>s$, and hence $C_1\cap C_2=\emptyset$. Finally, let
$e_L:=e_{L_n}(A,B_1\cup B_2)$. Then
\begin{equation}\label{eq:Ln-decomposition}
e(L_n)=e(L_n[A])+e_L+e(L_n[B_1\cup B_2]).
\end{equation}

We compare $L_n$ with $L'_n:=H_n$ from Lemma~\ref{lem:key_construction} to obtain a recursive bound for $\phi(n)$. Let $V(L'_n)=X_{n_0}\cup Y_{n_1}$ be its defining partition, and let
$B'_1\subseteq X_{n_0}$ and $B'_2\subseteq Y_{n_1}$ be the corresponding
sets from Lemma~\ref{lem:key_construction}. Then
$|B'_1|=|B'_2|=N$ and
$L'_n[B'_1,B'_2]\cong K_{N,N}$.
Set $A':=V(L'_n)\setminus(B'_1\cup B'_2)$, and let
$e_{L'}:=e_{L'_n}(A',B'_1\cup B'_2)$. 
Since all edges between $X_{n_0}$ and $Y_{n_1}$ are present and $e(L'_n[B'_1,X_{n_0} \setminus B'_1])=e(L'_n[B'_2,Y_{n_1} \setminus B'_2])=0$, we have
$e_{L'}
=(n-2N)N$. Hence 
\begin{equation}\label{eq:Lprime-decomposition}
e(L'_n)=e(L'_n[A'])+e_{L'}+e(L'_n[B'_1\cup B'_2])=f(n-2N,R)+(n-2N)N+e(L'_n[B'_1\cup B'_2]). \end{equation}
In the special case $R\cong2C_6$, we assert that $f(2N,2C_6)=g(N,N,2C_6)$. Let $q(t):=6\left\lfloor t/4\right\rfloor+\binom{t\bmod 4}{2}$. Since $q(t)\le {3t}/{2}$ for $t\ge 0$, 
for any $a+b=2N$,
$g(a,b,2C_6)\le ab+3N\le N^2+3N=g(N,N,2C_6)$, 
where the last equality uses $4\mid N$. Thus
$f(2N,2C_6)=g(N,N,2C_6)$, as asserted. 
By Lemma~\ref{lem:key_construction}, $e(L'_n[B'_1\cup B'_2]) +1 \ge g(N,N,R)$. 
Moreover, as shown above, $\Delta(L_n[B_i])\le s-1$ for $i=1,2$, and $L_n[B_1\cup B_2]$ is $\widehat{R}$-free. Therefore, applying the appropriate lemma among Lemmas~\ref{lem:special-two-C6},~\ref{lem:special-all-C4}, and~\ref{lem:bipartite}, with $(L_n[B_1\cup B_2],B_1,B_2)$ playing the role of $(G,V_0,V_1)$, we obtain
\begin{equation}\label{eq:local-bound}
e(L_n[B_1\cup B_2]) \le g(N,N,R) \le e(L'_n[B'_1\cup B'_2]) +1. \end{equation}
For sufficiently large $n$, we have $n-2N>M_0$. 
Since $L_n[A]$ is
$\widehat R$-free, we have
$e(L_n[A])\le \operatorname{ex}(n-2N,\widehat R)$.
Combining $\phi(n) =\operatorname{ex}(n,\widehat R)-f(n,R) =e(L_n)-e(L'_n)$, \eqref{eq:Ln-decomposition}, \eqref{eq:Lprime-decomposition}, and \eqref{eq:local-bound}, we
obtain
\begin{equation}\label{eq:phi-recursion}
\begin{aligned}[b]
\phi(n) & 
\le \operatorname{ex}(n-2N,\widehat R)-f(n-2N,R)+(e_L-e_{L'})+1 = \phi(n-2N)+(e_L-e_{L'})+1.
\end{aligned}
\end{equation}

If $e_{L'}-e_L\ge 2$, then \eqref{eq:phi-recursion} gives
$\phi(n)\le\phi(n-2N)-1$, so \ref{a} holds.
Suppose that some vertex $x\in V(L_n)$ satisfies
$d_{L_n}(x)<n/2+h(R)$. Since $L_n-x$ is $\widehat R$-free, we
have $e(L_n)-d_{L_n}(x)=e(L_n-x)\le \operatorname{ex}(n-1,\widehat R)$. Hence, by
 Lemma~\ref{lem:f-f},
$e(L_n)-\operatorname{ex}(n-1,\widehat R)\le d_{L_n}(x)<n/2+h(R)
\le f(n,R)-f(n-1,R)$. 
So $\phi(n)
<\phi(n-1)$. 
Since $\phi$ is integer-valued, \ref{b} holds.
We may therefore assume that neither \ref{a} nor \ref{b}
holds. Then $d_{L_n}(x)\ge n/2+h(R)$ for every
$x\in V(L_n)$ and $e_{L'}-e_L\le1$. It remains to prove~\ref{c}.

We now establish four structural estimates for the partition. 
Suppose that $\Delta(L_n[B_i\cup C_i])\ge s$ for some
$i\in\{1,2\}$. There exists a vertex $x\in B_i\cup C_i$ with $d_{L_n[B_i\cup C_i]}(x)\ge s$.
Choose an $s$-element set $ X\subseteq N_{L_n[B_i\cup C_i]}(x)$. By the definition of $C_i$ and the fact that $L_n[B_1,B_2]\cong K_{N,N}$, every vertex in $\left\{x\right\}\cup X$ has at most $\varepsilon N$ non-neighbors in $B_{3-i}$. Since $\varepsilon=2s/N$ and $N$ is sufficiently large,  $|N_{L_n[B_{3-i}]}(\left\{x\right\}\cup X)|\ge N-(s+1)\varepsilon N>s$. 
Therefore, we can construct a copy of $K_1+K_{s,s}$, which contains $\widehat{R}$, a contradiction. Thus, 
\begin{equation}\label{eq:BiCi-degree}
\Delta\bigl(L_n[B_i\cup C_i]\bigr)\le s-1
\qquad\text{for } i\in\{1,2\}.
\end{equation}
For every $x\in D$, there exists $i\in\{1,2\}$ such that $d_{L_n[B_i]}(x)<s$; otherwise, $x$ together with suitable $s$-subsets of $B_1$
and $B_2$ would contain a copy of $K_1+K_{s,s}$.
Since $x\notin C_i$, $d_{L_n[B_{3-i}]}(x)<(1-\varepsilon)N$, and hence $d_{L_n[B_1\cup B_2]}(x)<s+(1-\varepsilon)N = N-s$.
By \eqref{eq:BiCi-degree}, $e_{L_n}(B_i,C_i)\le N(s-1)$, and hence
$e_L \le 2N\cdot (s-1)+N\cdot(|C_1|+|C_2|)+|D|\cdot(N-s)
= e_{L'}+2N\cdot (s-1)-|D|\cdot s$.
Since $e_{L'}-e_L\le1$, $|D|\le {(1+2N(s-1))}/{s}$. Thus, setting $N_1:=\left\lceil(1+2N(s-1))/s\right\rceil$, we have
\begin{equation}\label{eq:D-bounded}
|D|\le N_1.
\end{equation}
By \eqref{eq:BiCi-degree}, for each $i\in\{1,2\}$, $e(L_n[B_i\cup C_i])\le \left\lfloor(s-1)|B_i\cup C_i|/2 \right\rfloor =O(n)$.
By \eqref{eq:D-bounded}, the number of edges incident with $D$ is $O(n)$.
Thus
$|B_1\cup C_1|\cdot|B_2\cup C_2|\ge e_{L_n}(B_1\cup C_1,B_2\cup C_2)\ge e(L_n)-O(n)\ge n^2/4-O(n)$.
Since $|B_1\cup C_1|+|B_2\cup C_2|=n-|D|$, 
$\bigl||B_1\cup C_1|-|B_2\cup C_2|\bigr|
=O(\sqrt n)$. Therefore, there exists a constant $N_2=N_2(s,N)$ such that
\begin{equation}\label{eq:Ci-size}
\bigl||B_i\cup C_i|-n/2\bigr|
\le N_2\sqrt n
\qquad\text{for }i\in\{1,2\}.
\end{equation}
Fix $i\in\{1,2\}$ and $x\in B_i\cup C_i$. By \eqref{eq:BiCi-degree}, $d_{L_n[B_i\cup C_i]}(x)\le s-1$. Since
$d_{L_n}(x)\ge n/2+ h(R)$, we have $d_{L_n[B_{3-i}\cup C_{3-i}\cup D]}(x) \ge n/2+h(R)-(s-1)$.
By \eqref{eq:D-bounded} and~\eqref{eq:Ci-size}, we have $|B_{3-i}\cup C_{3-i}\cup D|\le n/2+N_2\sqrt n+N_1$. 
Therefore, the number of non-neighbors of $x$ in $B_{3-i}\cup C_{3-i}\cup D$ is at most $N_2 \sqrt{n}+N_1+s-1-h(R)$. 
Thus, there exists a constant $N_3=N_3(R)$ such that 
\begin{equation}\label{eq:few-nonneighbors}
|B_{3-i}\cup C_{3-i}\cup D|
-d_{L_n[B_{3-i}\cup C_{3-i}\cup D]}(x)
\le N_3\sqrt n \qquad\text{for  } i\in\{1,2\} \text{ and } x\in B_i\cup C_i.
\end{equation}

We claim that $D$ admits a partition $D=D_1\cup D_2$ such that, for each
$i\in\{1,2\}$ and each $x\in D_i$, 
$d_{L_n[B_i\cup C_i]}(x)\le s-1$ and $d_{L_n[B_{3-i}\cup C_{3-i}]}(x)\ge n/6$.
Let $x\in D$. 
There exists $i(x)\in\left\{1,2\right\}$ such that $d_{L_n[B_{i(x)}\cup C_{i(x)}]}(x)\ge n/6$; otherwise, $d_{L_n}(x)<N_1+2\cdot ({n}/{6})<{n}/{2}+h(R)$, a contradiction. 
Suppose $x$ has at least $s$ neighbors in $B_{3-i(x)}\cup C_{3-i(x)}$. Choose a set $X$ of $s$ such neighbors.  By \eqref{eq:few-nonneighbors}, there are at least ${n}/{6}-sN_3\sqrt{n}>s$ common neighbors of $\left\{x\right\} \cup X$ in $B_{i(x)}\cup C_{i(x)}$. Hence $L_n$ contains a copy of $K_1+K_{s,s}$, which contains $\widehat{R}$, a contradiction. Therefore $d_{L_n[B_{3-i(x)}\cup C_{3-i(x)}]}(x)\le s-1$, so we put $x$ into $D_{3-i(x)}$.
This defines a partition $D=D_1\cup D_2$.

We now construct the partition needed for the key lemmas.
For $i\in\{1,2\}$, let $V_i:=B_i\cup C_i\cup D_i$.
By \eqref{eq:BiCi-degree} and the construction of $D_1,D_2$,
every vertex $x\in V_i$ satisfies
$d_{L_n[V_i]}(x)\le s-1+N_1$.
Hence $d_{L_n[V_{3-i}]}(x) \ge n/2+h(R)-(s-1)-N_1$.
By \eqref{eq:D-bounded} and \eqref{eq:Ci-size}, there
exists a constant $N_4=N_4(s,N)$ such that
$\left||V_i|- n/2\right|\le N_4\sqrt n$ for each $i\in\{1,2\}$.
Thus, for every $x\in V_i$,
$|V_{3-i}|-d_{L_n[V_{3-i}]}(x)
\le N_4\sqrt n+s-1+N_1-h(R)$.
Suppose that some $x\in V_i$ has $s$ distinct neighbors in $V_i$; let $X'$ denote this set. Since
$|B_{3-i}\cup C_{3-i}|
\ge
 n/2-N_2\sqrt n$, 
$|N_{L_n[B_{3-i}\cup C_{3-i}]}(X'\cup \{x\})|\ge n/2-N_2\sqrt n
-(s+1)\bigl(N_4\sqrt n+s-1+N_1-h(R)\bigr)\ge s$.
Thus $L_n$ contains a copy of $K_1+K_{s,s}$ and hence a copy
of $\widehat{R}$, a contradiction. Thus $V(L_n)=V_1\cup V_2$  is a partition satisfying $\Delta(L_n[V_i])\le s-1$ and $|V_i|=n/2+O(\sqrt{n})$ for $i\in \{1,2\}$. 
Relabel $V_1$ and $V_2$, if necessary, so that $|V_1| \ge |V_2|$. Depending on the form of $R$, we apply the appropriate lemma among Lemmas~\ref{lem:special-two-C6}, \ref{lem:special-all-C4},
 and~\ref{lem:bipartite}. Together with $e(L_n)\ge f(n,R)$, this yields $e(L_n)=f(n,R)$ and $L_n\in\mathcal K_n(R)$,  hence \ref{c} holds.
\end{proof}

\section{Proofs of the key lemmas}\label{sec:proof-lemmas}

\subsection{Auxiliary lemmas}

We first record several auxiliary lemmas used in the proofs of the key lemmas.

\begin{lemma}[McKay--Wormald--Wysocka~\cite{MWW2004}]\label{lem:mww-short-cycles}
For every fixed integer $d\ge3$ and every fixed finite set of cycle lengths, for all sufficiently large $n$, if $nd$ is even, then there exists a $d$-regular graph on $n$ vertices containing no cycle of any of these lengths.
\end{lemma}

\begin{lemma}[Yuan~\cite{Yuan2021}]\label{lem:many-paths,stars}
For every pair of integers $k,k'$ with $k>k'\ge0$, every positive integer $q$, and every $\varepsilon>0$, there exists a constant $n_0=n_0(k,k',q,\varepsilon)$ such that the following holds. Let $H$ be an $n$-vertex graph with $n\ge n_0$, $\Delta(H)\le k-1$, and $e(H)\ge(k'n+\varepsilon n)/2$. Then $H$ contains $q\cdot P_{k'+2}$ and $q\cdot S_{k'+2}$ as subgraphs.
\end{lemma}

\begin{lemma}[Yuan~\cite{Yuan2021}]
\label{lem:yuan-star-path}
Let $k\ge2$ and $n\ge2k$. Then $\operatorname{ex}\bigl(n,\left\{S_{k+1},P_{2k-1}\right\}\bigr)
=
\left\lfloor{(k-1)n}/{2}\right\rfloor$.
\end{lemma}

\begin{lemma}[Faudree--Schelp~\cite{FaudreeSchelp1975},
Kopylov~\cite{Kopylov1977}]\label{lem:ex(P_5)}
For every $n$,
$\operatorname{ex}(n,P_5)
=
6\left\lfloor n/4\right\rfloor
+\binom{n\bmod4}{2}$.
Moreover, every extremal $P_5$-free graph is isomorphic to
$\left\lfloor n/4\right\rfloor K_4\cup K_{n\bmod4}$.
\end{lemma}

\begin{lemma}[Chv\'atal--Hanson~\cite{ChvatalHanson1976}]\label{lem:chvatal-hanson}
Let $\nu(G)$ denote the matching number of $G$.
For integers $\nu,\Delta\ge1$, if a graph $H$ satisfies
$\nu(H)\le\nu$ and $\Delta(H)\le\Delta$, then
$e(H)
\le \nu\Delta+\nu$.
\end{lemma}

\begin{lemma}\label{lem:path-cycle}
Let $p\ge2$ and $1\le k\le p$ be integers. Let $G$ be a graph with a partition $V(G)=V_0\cup V_1$, where $|V_0|= k$ and $G[V_0,V_1]\cong K_{k,|V_1|}$. If $G[V_1]$ contains $k$ pairwise vertex-disjoint paths whose total number of vertices is at least $2p-k$, then $G$ contains a copy of $C_{2p}$.
In particular, \begin{enumerate}[label=\textup{(\alph*)}]
\item\label{lem:path-cycle-p3}
if $\lceil p/2\rceil\le c\le p$, then
$cK_1+\bigl((p-c)P_3\cup(2c-p)P_1\bigr)$ contains a copy of $C_{2p}$;
\item\label{lem:path-cycle-p2}
if $p\ge3$, then
$(p-1)K_1+\bigl(2P_2\cup(p-3)P_1\bigr)$ contains a copy of $C_{2p}$.
\end{enumerate}
\end{lemma}

\begin{proof}
By truncating the paths in $G[V_1]$ if necessary, we obtain $k$ paths with
exactly $2p-k$ vertices in total. As $2p-k\ge k$ and $G[V_0,V_1]$ is complete, these paths together with $V_0$ form a copy of $C_{2p}$. For~\ref{lem:path-cycle-p3} and~\ref{lem:path-cycle-p2}, apply
the main statement with $k=c$ and $k=p-1$, respectively;
in the two cases the corresponding path systems have
$2p-c$ and $p+1=2p-(p-1)$ vertices in total.
\end{proof}

The following two elementary integer lemmas are needed in the proof of Lemma~\ref{lem:bipartite}.

\begin{lemma}\label{lem:int1}
Let $\mathbf p=(p_1,\ldots,p_r)$ be a nondecreasing sequence of
integers with $p_1\ge2$, and set $s:=\sum_{i=1}^r p_i$. For every integer $k$ with
$2\le k\le s-3$, unless
$(k,\mathbf p)=(3,(3,3))$, at least one of the following holds: 
\begin{enumerate*}[label=\textup{(\arabic*)}]
\item\label{first_desired} 
$\sum_{i=1}^r\lceil 2p_i/(s+1-k)\rceil\le k$, 
\item\label{second_desired} $\sum_{i=1}^r\lceil 2p_i/(k+2)\rceil\le s-k$.
\end{enumerate*}
\end{lemma}

\begin{proof}
Let $a=s+1-k$ and $b=k+2$. Since $2\le k\le s-3$, we have $a\ge4$, $b\ge4$, and $a+b=s+3$.

\noindent\textbf{Case 1.} $r\ge3$. 

For each $i\in[r]$, since $r\ge3$ and $p_j\ge2$ for every
$j\in[r]$, we have $a+b=3+ p_i+\sum_{j\neq i}p_j\ge p_i+7$.
Suppose that both \ref{first_desired} and \ref{second_desired} failed. Since the two sums are integers, we would have $\sum_{i=1}^r\left(\left\lceil 2p_i/a\right\rceil+\left\lceil 2p_i/b\right\rceil\right)\ge s+2$. So it suffices to prove that $\sum_{i=1}^r\left(\left\lceil 2p_i/a\right\rceil+\left\lceil 2p_i/b\right\rceil\right)\le s+1$.
We claim that 
if $p\ge2$, $a,b\ge4$, and $a+b\ge p+7$, then $\left\lceil 2p/a\right\rceil+\left\lceil 2p/b\right\rceil\le p$, except when $(a,b,p)=(5,5,3)$, in which case $\left\lceil 2p/a\right\rceil+\left\lceil 2p/b\right\rceil\le p+1$.
Indeed, assume without loss of generality that $a\le b$. Then $b\ge(p+7)/2$. 
If $p\ge5$, then $b\ge6$, and hence $\left\lceil 2p/a\right\rceil+\left\lceil 
2p/b\right\rceil\le 
\left\lceil p/2\right\rceil+\left\lceil p/3\right\rceil \le p$. 
If $p\in \left\{2,4\right\}$, then $\left\lceil 2p/a\right\rceil+\left\lceil 2p/b\right\rceil\le p$.
It remains to consider $p=3$. In this case, $a+b\ge10$.
 If $(a,b)\neq(5,5)$, then $b \ge 6$, so $\left\lceil 2p/a\right\rceil+\left\lceil 2p/b\right\rceil\le p$. If $(a,b)=(5,5)$, then $\left\lceil 2p/a\right\rceil+\left\lceil 2p/b\right\rceil=p+1$, as claimed. 
When $(a,b)=(5,5)$, since $r\ge 3$ and $s=7$, we have $\mathbf p=(2,2,3)$ and so $\sum_{i=1}^r\left(\left\lceil 2p_i/a\right\rceil+\left\lceil 2p_i/b\right\rceil\right)\le s+1$. 
When $(a,b)\ne (5,5)$, we have $\sum_{i=1}^r\left(\left\lceil 2p_i/a\right\rceil+\left\lceil 2p_i/b\right\rceil\right)\le \sum_{i=1}^r p_i=s$, as desired.

\noindent\textbf{Case 2.} $r=2$, $(k,\mathbf{p}) \neq (3,(3,3))$. 

If $6\le k\le s-3$, then $s\ge9$. Using $\lceil x/y\rceil=\left\lfloor(x-1)/y\right\rfloor+1$ for positive integers $x,y$, we obtain $\left\lceil 2p_1/a\right\rceil+\left\lceil 2p_2/a\right\rceil\le(2s-2)/(s+1-k)+2\le k$, so \ref{first_desired} holds.
When $k\in\left\{4,5\right\}$, recall that $s=a+k-1$. If $p_2\le a$, then $\left\lceil 2p_1/a\right\rceil+\left\lceil 2p_2/a\right\rceil\le4\le k$. If $p_2\ge a+1$, then $p_1=s-p_2\le k-2$. Since $p_2=s-p_1\le a+k-3$ and $a\ge 4$, we have 
$\left\lceil 2p_1/a\right\rceil+\left\lceil 2p_2/a\right\rceil
\le\left\lceil 2(k-2)/a\right\rceil+\left\lceil2(a+k-3)/a\right\rceil
\le k$, and hence \ref{first_desired} holds.
When $k=3$, we have $s\ge6$. Since
$(p_1,p_2)\ne(3,3)$ and $2\le p_1\le p_2$, we have $p_2\ge4$. For every integer $p\ge2$, we have $\left\lceil 2p/5\right\rceil\le p-1$, while for every integer $p\ge4$, we have $\left\lceil 2p/5\right\rceil\le p-2$. Hence $\left\lceil 2p_1/5\right\rceil+\left\lceil 2p_2/5\right\rceil\le p_1+p_2-3=s-k$, which proves \ref{second_desired}. 
When $k=2$, we have $\left\lceil 2p_1/(k+2)\right\rceil+\left\lceil 2p_2/(k+2)\right\rceil\le(p_1-1)+(p_2-1)=s-k$, so \ref{second_desired}  holds. 

\noindent\textbf{Case 3.} $r=1$. 

If $3\le k\le s-3$, then $s=p_1\ge6$ and $k(p_1+1-k)\ge2p_1$. Since $p_1+1-k> 0$ and $k \in \mathbb{Z}$, we have $k \ge \left\lceil 2p_1/(s+1-k)\right\rceil$, so \ref{first_desired} holds.
If $k=2$, then $p_1\ge5$, and $s-k=p_1-2\ge p_1/2=2p_1/(k+2)$. Since $s-k \in \mathbb{Z}$, we have $s-k\ge\left\lceil 2p_1/(k+2)\right\rceil$, and hence \ref{second_desired} holds.
\end{proof}

\begin{lemma}\label{lem:int2}
Let $\mathbf{p}=(p_1,\ldots,p_r)$ be a nondecreasing sequence of
integers with $p_1\ge2$, and set $s:=\sum_{i=1}^r p_i$. Let $k$ be an integer satisfying
$2\le k\le s-3$. Suppose that 
\begin{enumerate*}[label=\textup{(\roman*)}]
\item\label{f} 
$\sum_{i=1}^r\lceil 2p_i/(k+2)\rceil\le s-k$,
\item\label{g} $\sum_{i=1}^r\lceil 2p_i/(s+1-k)\rceil>k$,
\item\label{h} 
$\sum_{i=1}^r\left\lfloor p_i/2\right\rfloor<k$.
\end{enumerate*}
Then $(k,\mathbf{p})\in
\left\{(4,(2,3,3)),(5,(3,3,3))\right\}$.
\end{lemma}

\begin{proof}
Both listed choices of $(k,\mathbf p)$ clearly satisfy \ref{f}, \ref{g} and \ref{h}. Since $p_i\ge2$ for every $i\in[r]$, we have
$r\le\sum_{i=1}^r\lfloor p_i/2\rfloor<k$.
Set $c:=k-r$. Then $c\ge1$.
By \ref{f}, we have $r\le s-k$ and hence $s\ge 2r+c$. Set $a:=s+1-k=s+1-r-c$. Since $k\le s-3$, we have $a\ge4$.
By \ref{g}, we have $k+1\le\sum_{i=1}^r\lceil 2p_i/a\rceil$. Using
$\lceil m/a\rceil\le(m-1)/a+1$ for positive integers $m$ and $a$, we obtain
$k+1\le r+(2s-r)/a$, and therefore
$c+1\le(2s-r)/(s+1-r-c)$. Rearranging yields
\begin{equation}\label{eq:c-s-bound}
(c-1)s\le cr+c^2-1.
\end{equation}
We next show that $r\ge2$. Suppose that $r=1$. Then $s=p_1$. If
$k\ge3$, then the argument in the $r=1$ case of Lemma~\ref{lem:int1} gives $\lceil 2p_1/(p_1+1-k)\rceil\le k$, contradicting \ref{g}. If $k=2$, then $\lfloor p_1/2\rfloor<2$ and hence $p_1\le3$,
contradicting $k\le s-3=p_1-3$. Hence $r\ge2$.

\noindent\textbf{Case 1.} $c \ge 3$. 

Combining $s\ge2r+c$ with \eqref{eq:c-s-bound}, we obtain $(c-2)r\le c-1$.
If $c\ge4$, then $r\le3/2$, contradicting $r\ge2$. If $c=3$, then $r\le2$ and hence $k=5$. 
However, \eqref{eq:c-s-bound} yields $s \le 7$, which contradicts $k\le s-3$.

\noindent\textbf{Case 2.} $c=2$.

In this case, $k=r+2$. 
From \eqref{eq:c-s-bound} and $s\ge2r+c$, we obtain
$2r+2\le s\le2r+3$.
Suppose first that $s=2r+2$. Then
$\mathbf p\in\{(2,\ldots,2,4),(2,\ldots,2,3,3)\}$.
Since $k\le s-3$ and $s-k=r$, we have $r\ge3$. Thus $a=r+1\ge4$ and hence $\sum_{i=1}^r\lceil 2p_i/a\rceil\le r+2=k$.
Suppose now that $s=2r+3$. Then $\mathbf{p}\in \{(2,\ldots,2,5), (2,\ldots,2,3,4), (2,\ldots,2,3,3,3)\}$. 
Since $r\ge2$, we have $a=r+2\ge4$.
For $\mathbf{p}\in \{(2,\ldots,2,5),(2,\ldots,2,3,4)\}$, we have $\sum_{i=1}^r\lceil 2p_i/a\rceil\le r+2=k$.
For $\mathbf p=(2,\ldots,2,3,3,3)$, we have $r\ge3$. If $r\ge4$, then $a=r+2\ge6$. Hence $\sum_{i=1}^r\lceil 2p_i/a\rceil=r<k$. All cases considered above contradict \ref{g}, except when
$\mathbf p=(3,3,3)$, which gives $k=5$.

\noindent\textbf{Case 3.} $c=1$. 

In this case, $k=r+1$. 
Since $\sum_{i=1}^r\lfloor p_i/2\rfloor<r+1$ and each $p_i\ge2$, we must have $p_i\in\{2,3\}$ for every $i$. Let $n_j:=|\{i:p_i=j\}|$ for $j\in\{2,3\}$. Then $r=n_2+n_3$, $s=2n_2+3n_3$, and $a=r+n_3$.
If $a\ge6$, then $\sum_{i=1}^r\lceil 2p_i/a\rceil=r<k$, contradicting \ref{g}.
Thus $a\le5$ and hence $r\in\{2,3,4,5\}$.
Since $k\le s-3$, $a\ge 4$. Note that $a=r+n_3$ and $n_3\le r$. 
If $r=2$, then $n_3=2$ and hence $\sum_{i=1}^r\lceil 2p_i/(k+2)\rceil=4>s-k$, contradicting \ref{f}. If $r\in \left\{4,5\right\}$, then $n_3\le1$ and $4\le a\le 5$. Hence $\sum_{i=1}^r\lceil 2p_i/a\rceil=n_2+2n_3 \le k$, contradicting \ref{g}. 
Now consider $r=3$. 
Since $4\le r+n_3\le5$, $n_3\in\{1,2\}$. 
If $(n_2,n_3)=(2,1)$, then $\sum_{i=1}^r\lceil 2p_i/a\rceil=4=k$, contradicting \ref{g}. Hence $(n_2,n_3)=(1,2)$, giving $(k,\mathbf p)=(4,(2,3,3))$.
\end{proof}

The next lemma is a local embedding lemma. It will be used in the proof of Lemma~\ref{lem:bipartite}.

\begin{lemma}\label{lem:local-embedding-from-C1}
Let $R\cong \bigcup_{i=1}^r C_{2p_i}$ with $2\le p_1\le\cdots\le p_{r}$ and let $s=\sum_{i=1}^r p_i$. Let $G$ be a graph with a partition $V(G)=V_0\cup V_1$, where $|V_0|=n_0\ge n_1=|V_1|\ge s+2$. Set $G_i:=G[V_i]$ and $k_i:=\Delta(G_i)\le s-1$ for
$i\in\{0,1\}$. For every fixed positive integer $N$, if $n_0$ is sufficiently large and
\[e(G_0)\le \left\lfloor\frac{(s-1-k_1)n_0}{2}\right\rfloor+N,\quad
e(G)\ge n_0n_1+\left\lfloor\frac{(s-1)n_0}{2}\right\rfloor+1, \quad
\sum_{i=1}^r\left\lceil\frac{2p_i}{k_1+2}\right\rceil\le s-k_1,\]
then $G$ contains a copy of $\widehat{R}$.
\end{lemma}

\begin{proof}
We split into two cases according to whether $n_1\ge n_0/2$.

\noindent\textbf{Case 1.} $n_1\ge n_0/2$.

We assert that there is a vertex $y\in V_0$ with $d_{G_0}(y)\ge s-k_1$.
Since $\Delta(G_1)=k_1$ and $n_1\le n_0$, we have
$e(G_1)\le\lfloor k_1n_1/2\rfloor\le\lfloor k_1n_0/2\rfloor$.
Using $\lfloor a\rfloor-\lfloor b\rfloor\ge\lfloor a-b\rfloor$, we obtain
$e(G_0)
\ge e(G)-e(G_1)-n_0n_1
\ge
\left\lfloor{(s-1-k_1)n_0}/{2}\right\rfloor+1$. 
Thus $\Delta(G_0)\ge s-k_1$, and such a vertex $y$ exists.
On the other hand, $\Delta(G_0)\le s-1$ gives
$e(G_0)\le\lfloor(s-1)n_0/2\rfloor$. Hence $k_1\ge1$.
Note that $\lceil 2p_i/(k_1+2)\rceil\le p_i$ for every $i\in[r]$ and
$\sum_{i=1}^r\lceil 2p_i/(k_1+2)\rceil\le s-k_1\le\sum_{i=1}^r p_i$.
Hence we can choose integers $c_1,\ldots,c_r$ such that
$\lceil 2p_i/(k_1+2)\rceil\le c_i\le p_i$ for every $i\in[r]$ and
$\sum_{i=1}^r c_i=s-k_1$. 
We next find suitable paths in $G_1$ for embedding $R$.
Using the assumed upper bound on $e(G_0)$ and $n_1\le n_0$, we obtain 
\[
e(G_1)
\ge e(G)-e(G_0)-n_0n_1
\ge
\left\lfloor\frac{(s-1)n_0}{2}\right\rfloor
-\left\lfloor\frac{(s-1-k_1)n_0}{2}\right\rfloor
-N+1
\ge
\left\lfloor\frac{k_1n_1}{2}\right\rfloor-N+1.
\]
Since $n_1\ge n_0/2$ and $n_0$ is sufficiently large, we have
$e(G_1)>(k_1-1/2)n_1/2$. Since $k_1\ge1$, Lemma~\ref{lem:many-paths,stars} with $k=s$, $k'=k_1-1$, and
$\varepsilon=1/2$ yields a copy of $N(s-k_1)\cdot P_{k_1+1}$ in $G_1$.
Partition these paths into $N$ vertex-disjoint subgraphs $Q_1,\ldots,Q_N$, each consisting of $s-k_1$ copies of
$P_{k_1+1}$. Since $\sum_{j=1}^r c_j=s-k_1$, for each $i\in[N]$ we
may decompose
$Q_i=Q_{i1}\mathbin{\cup}\cdots\mathbin{\cup}Q_{ir}$ so that
$Q_{ij}\cong c_j\cdot P_{k_1+1}$ for every $j\in[r]$.

Suppose first that all edges between $V(Q_\ell)$ and $N_{G_0}[y]$
are present for some $\ell\in[N]$. Since
$\sum_{i=1}^r c_i=s-k_1\le d_{G_0}(y)$, we may choose pairwise
disjoint sets $Z_1,\ldots,Z_r\subseteq N_{G_0}(y)$ with
$|Z_i|=c_i$ for every $i\in[r]$.  Moreover, by the choice of $c_i$, $|V(Q_{\ell i})|=c_i(k_1+1)\ge2p_i-c_i$.
By Lemma~\ref{lem:path-cycle},
$G[Z_i\cup V(Q_{\ell i})]$ contains a copy of $C_{2p_i}$ for each $i$.
These cycles are pairwise vertex-disjoint, and since $y$ is adjacent to all their vertices, their union together with $y$ contains a copy of $\widehat R$.
If no such $\ell$ exists, then for every $\ell\in[N]$, there is a
missing edge between $V(Q_\ell)$ and $N_{G_0}[y]$. Since the sets $V(Q_1),\ldots,V(Q_N)$ are pairwise disjoint, these $N$ missing edges are distinct. Hence
$e_G(V_0,V_1)\le n_0n_1-N$. Consequently,
$e(G)
\le
\left\lfloor{(s-1-k_1)n_0}/{2}\right\rfloor+N
+\left\lfloor{k_1n_1}/{2}\right\rfloor
+n_0n_1-N
<
n_0n_1+\left\lfloor{(s-1)n_0}/{2}\right\rfloor+1$, a contradiction. Hence $G$ contains a copy of $\widehat R$.

\noindent\textbf{Case 2.} $n_1<n_0/2$.

In this case, $e(G_1)\le\lfloor k_1n_1/2\rfloor
\le\lfloor k_1n_0/4\rfloor$. Hence $e(G)
\le
\left\lfloor{(s-1-k_1)n_0}/{2}\right\rfloor+N
+
\left\lfloor {k_1n_0}/{4} \right\rfloor
+n_0n_1$. Note that $N$ is fixed and $n_0$ is sufficiently
large. 
If $k_1\ge1$, then
$e(G)
< n_0n_1+\left\lfloor{(s-1)n_0}/{2}\right\rfloor+1$.
If $k_1=0$, then $e(G_1)=0$ and 
$e(G)\le \left\lfloor{(s-1)n_0}/{2}\right\rfloor+n_0n_1
< \left\lfloor{(s-1)n_0}/{2}\right\rfloor+n_0n_1 +1$.
In either case, this contradicts the assumed lower bound on $e(G)$.
\end{proof}

\subsection{The proofs of construction lemmas}

\begin{proof}[\textbf{Proof of Lemma~\ref{lem:U-nonempty}}]
We distinguish two cases according to $r$.

\noindent\textbf{Case 1.} $r=1$. 

In this case, $R\cong C_{2p_1}$, $s=p_1$, and
$R^-=\{P_{2p_1-1}\}$. By Lemma~\ref{lem:yuan-star-path},
for every sufficiently large $n$,
$\operatorname{ex}\bigl(n,\{S_{p_1+1},P_{2p_1-1}\}\bigr)
=
\left\lfloor{(s-1)n}/{2}\right\rfloor$.
Hence there exists an $n$-vertex $P_{2p_1-1}$-free graph $H$ with
$\Delta(H)\le s-1$ and
$e(H)=\left\lfloor{(s-1)n}/{2}\right\rfloor$.
It follows that $H$ is nearly $(s-1)$-regular.

\noindent\textbf{Case 2.} $r\ge2$. 

Since every member of $R^-$ contains a cycle $C_{2p_i}$ for some $i\in[r]$, it suffices to construct an $n$-vertex nearly $(s-1)$-regular graph which is $\{C_{2p_1},\ldots,C_{2p_r}\}$-free. Note that $s-1\ge3$.
Suppose first that $(s-1)n$ is even. By Lemma~\ref{lem:mww-short-cycles}, for every sufficiently large $n$, there exists such an
$(s-1)$-regular graph $H$. 
Now suppose that $(s-1)n$ is odd. Then both $s-1$ and $n$ are odd, so $(s-1)(n+1)$ is even. By Lemma~\ref{lem:mww-short-cycles}, for every sufficiently large $n$ there exists an $(s-1)$-regular graph $H'$ on $n+1$ vertices containing no cycle of length at most $2s$. Fix $x\in V(H')$. Since $H'$ is triangle-free, $N_{H'}(x)$ is independent. Pair all but one of its $s-1$ vertices, and let $M$ be the resulting matching. The graph
$H:=(H'-x)\cup M$
is an $n$-vertex nearly $(s-1)$-regular graph.
We claim that $H$ contains none of the cycles $C_{2p_i}$. Suppose
otherwise, and let $C\cong C_{2p_i}$ be a cycle in $H$. Since
$p_i<s$ and $H'-x$ contains no cycle of length at most $2s$, the
cycle $C$ must contain an edge of $M$. Deleting all edges of
$M\cap E(C)$ from $C$ yields nontrivial paths in $H'-x$, since $M$
is a matching. Let $Q$ be one such path, with endpoints $u,v$. Then
$u,v\in N_{H'}(x)$, so $Q\cup\{ux,xv\}$ is a cycle in $H'$ of length
at most
$e(Q)+2\le 2p_i+1\le 2s$,
a contradiction. Hence $H$ is the desired graph.
\end{proof}

\begin{proof}[\textbf{Proof of Lemma~\ref{lem:construction-free}}]
Let $G\in\mathcal K_{n_0,n_1}(R)$. Write $G=G_0+G_1$ as in the definition of $\mathcal K_{n_0,n_1}(R)$, where $|V(G_i)|=n_i$ for $i\in\{0,1\}$. Then $\Delta(G_i)\le s-1$ for $i\in\{0,1\}$, and in the special case $R\cong2C_6$ we have the stronger bound
$\Delta(G_i)\le3$. Set $V_i:=V(G_i)$.
Suppose, for a contradiction, that $G$ contains a copy of $\widehat{R}$ with apex $c$ and cycle components
$Q_1,\ldots,Q_r$, where $Q_i\cong C_{2p_i}$ for $i\in[r]$.
Then $\bigcup_{i=1}^rV(Q_i)\subseteq N_G(c)$. Let $h\in\{0,1\}$ be such that $c\in V_h$.

\begin{claim}\label{Q_j}
There exists $j\in[r]$ such that $|V(Q_j) \cap V_h|\le p_j-1$. Moreover, when $R\cong2C_6$, there exists $j\in[2]$ such that $|V(Q_j) \cap V_h|\le 1$.
\end{claim}
\begin{proof}
Suppose, to the contrary, that $|V(Q_i)\cap V_h|\ge p_i$ for every $i\in[r]$. Since the cycles $Q_1,\ldots,Q_r$ are
pairwise vertex-disjoint and all their vertices are adjacent to
$c$, we have $d_{G_h}(c)\ge\sum_{i=1}^r |V(Q_i) \cap V_h|\ge \sum_{i=1}^r p_i=s$, contradicting $d_{G_h}(c)\le s-1$.
When $R\cong2C_6$, if $|V(Q_i) \cap V_h|\ge 2$ for both $i\in[2]$, then $d_{G_h}(c)\ge\sum_{i=1}^2 |V(Q_i) \cap V_h|\ge 4$, contradicting $d_{G_h}(c)\le 3$.
\end{proof}
\noindent\textbf{Case 1.} $R\cong2C_6$.

For each $\ell\in\{0,1\}$,
$G_\ell\cong
(\left\lfloor\frac{n_\ell}{4}\right\rfloor K_4
\cup K_{n_\ell\bmod4})$,
so $G_\ell$ is both $C_6$-free and $P_5$-free. Hence each $Q_i$ meets
both parts. By Claim~\ref{Q_j}, some $Q_j$ satisfies
$|V(Q_j)\cap V_h|=1$. Deleting this vertex from $Q_j$ leaves a copy
of $P_5$ in $G_{1-h}$, a contradiction.

\noindent\textbf{Case 2.} $R\not\cong2C_6$.

In this case, $G_0$ is nearly $(s-1)$-regular and $R^-$-free,
while $\Delta(G_1)\le1$. In particular, $G_1$ is acyclic, so every
$Q_i$ meets $V_0$.
Suppose first that $c\in V_1$. 
Since
$V_1\cap\bigcup_{i=1}^rV(Q_i)\subseteq N_{G_1}(c)$,
we have
$\left|V_1\cap\bigcup_{i=1}^rV(Q_i)\right|
\le d_{G_1}(c)\le1$.
Deleting its unique vertex if $V_1\cap\bigcup_{i=1}^rV(Q_i)\ne\emptyset$, and otherwise deleting an arbitrary vertex from one of the cycles, leaves a member of $R^-$ contained in $G_0$, a
contradiction.
It remains to consider $c\in V_0$. By Claim~\ref{Q_j}, there exists
$j\in[r]$ such that 
$|V(Q_j)\cap V_0|\le p_j-1$.
Set
$k':=|V(Q_j)\cap V_0|$.
Since $Q_j$ meets $V_0$, we have $1\le k'\le p_j-1$. 
If $R\cong rC_4$, then $p_j=2$, and hence $k'=1$. Deleting the
unique vertex of $V(Q_j)\cap V_0$ from $Q_j$ leaves a copy of $P_3$
in $G_1$, contradicting $\Delta(G_1)\le1$.
If $R\not\cong rC_4$, then $e(G_1)\le1$. Let
$F:=Q_j\bigl[V(Q_j)\cap V_1\bigr]$.
The graph $F$ is a linear forest on $2p_j-k'$ vertices with at most
$k'$ components. Therefore,
$e(G_1)\ge e(F)
\ge (2p_j-k')-k'
\ge2$,
contradicting $e(G_1)\le1$.

In all cases we get a contradiction. Therefore $G$ is $\widehat{R}$-free. 
Since graphs in $\mathcal K_n(R)$ have exactly $f(n,R)$ edges, the lower bound $\operatorname{ex}(n,\widehat{R})\ge f(n,R)$ follows for all sufficiently large $n$.
\end{proof}

\begin{proof}[\textbf{Proof of Lemma~\ref{lem:key_construction}}]
For any $N_0$,
choose $N=N(R,N_0)>N_0$ divisible by $4$,
subject to the additional divisibility conditions specified below. Let
$(n_0-N,n_1-N)$ satisfy $g(n_0-N,n_1-N,R)=f(n-2N,R)$. For sufficiently large $n$, both $n_0-N$ and $n_1-N$ are sufficiently large. By the definition of $\mathcal{K}_{n-2N}(R)$ and using 
Lemma~\ref{lem:U-nonempty} when necessary, choose
$H_{n-2N}\in\mathcal K_{n-2N}(R)$ with defining partition
$V(H_{n-2N})=X_{n_0-N}\cup Y_{n_1-N}$ as specified below. 
First, we prove the following claim.

\begin{claim}\label{claim:g-g}
Fix a positive integer $N$ and a bipartite $2$-regular graph $R$. Let   $n$ be sufficiently large.  Suppose that  there exists a quantity $C_N(n)$, independent of $a$ and $b$, such that
$
g(a,b,R)-g(a-N,b-N,R)=C_N(n)
$ for all integers $a,b\ge N+2$ with $a+b=n$.
Then 
$
g(a,b,R)=f(n,R)
$ if and only if $
g(a-N,b-N,R)=f(n-2N,R).
$
\end{claim}
\begin{proof}
For sufficiently large $n$,  
$f(n,R)=\max_{a+b=n} g(a,b,R)=\max_{a,b\ge N+2, a+b=n} g(a,b,R)
=C_N(n)+\max_{a,b\ge N+2, a+b=n} g(a-N,b-N,R)
=C_N(n)+f(n-2N,R)$.
Hence, for any $a,b\ge N+2$ with $a+b=n$,
$g(a,b,R)=f(n,R)$ if and only if
$g(a-N,b-N,R)+C_N(n)=f(n-2N,R)+C_N(n)$, which is equivalent to
$g(a-N,b-N,R)=f(n-2N,R)$.
\end{proof}

\noindent\textbf{Case 1.} $R\cong2C_6$.

Choose $H_{n-2N}$ so that $H_{n-2N}[X_{n_0-N}]\cong \lfloor (n_0-N)/4\rfloor K_4\cup K_{(n_0-N)\bmod4}$ and $H_{n-2N}[Y_{n_1-N}]\cong \lfloor (n_1-N)/4\rfloor K_4\cup K_{(n_1-N)\bmod4}$. 
Let $G_{R,0}\cong G_{R,1}\cong(N/4)K_4$ and $H'_n:=\left(H_{n-2N}[X_{n_0-N}]\cup G_{R,0}\right) + \left(H_{n-2N}[Y_{n_1-N}]\cup G_{R,1}\right)$. Set
$B'_1:=V(G_{R,0})$, $B'_2:=V(G_{R,1})$, $X_{n_0}:=X_{n_0-N}\cup B'_1$ and $Y_{n_1}:=Y_{n_1-N}\cup B'_2$. Then \ref{1},~\ref{2},~\ref{3} hold and $H'_n\cong \left(\lfloor n_0/4\rfloor K_4\cup K_{n_0\bmod4}\right)+\left(\lfloor n_1/4\rfloor K_4\cup K_{n_1\bmod4}\right)$. 
So $H'_n\in \mathcal{K}_{n_0,n_1} (R)$. 
Moreover, $\max\{\Delta(H'_n[X_{n_0}]),\Delta(H'_n[Y_{n_1}])\}
\le s-1$ and
$e(H'_n[B'_1\cup B'_2])=N^2+3N=g(N,N,R)$. Thus \ref{4} holds.
Since $4\mid N$, for all $a,b\ge N+2$ with $a+b=n$, 
$g(a,b,R)-g(a-N,b-N,R)
=N(n+3-N)$, which is independent of $(a,b)$.  
As $g(n_0-N,n_1-N,R)=f(n-2N,R)$,
Claim~\ref{claim:g-g} gives $g(n_0,n_1,R)=f(n,R)$.  Therefore, $e(H'_n)=f(n,R)$ and $H'_n\in\mathcal K_n(R)$. 

\noindent\textbf{Case 2.} $R\not\cong 2C_6$.

Choose $H_{n-2N}$ so that $H_{n-2N}[X_{n_0-N}]$ is nearly $(s-1)$-regular and $\{C_{2p_1},\ldots,C_{2p_r}\}$-free; when $r=1$, it is also
$P_{2p_1-1}$-free. Moreover,  $H_{n-2N}[Y_{n_1-N}]$ is nearly $1$-regular when $R\cong rC_4$, while 
it has exactly one edge when $R\not\cong rC_4$.
Let $G_{R,2}:=(N/2)K_2$ when $R\cong rC_4$, while $G_{R,2}:=NK_1$ when $R\not\cong rC_4$. 
For $r=1$, 
we additionally require that $4p_1\mid N$ and
let $G_{R,1}:=(N/p_1)K_{p_1}$. Then $G_{R,1}$ is
$(p_1-1)$-regular and $P_{2p_1-1}$-free.
For $r\ge2$, by Theorem~\ref{thm:erdos-sachs}, there
exists an $(s-1)$-regular graph on some
$N'\le4\sum_{i=1}^{2p_r-1}(s-2)^i$ vertices with girth at least
$2p_r+1$. We additionally require that $4N'\mid N$, and let $G_{R,1}$ be the disjoint union of $N/N'$ copies of this graph. Then
$G_{R,1}$ is $(s-1)$-regular and $\{C_{2p_1},\ldots,C_{2p_r}\}$-free. 
Define $B_i':=V(G_{R,i})$ for $i\in\{1,2\}$.
Set $X_{n_0}:=X_{n_0-N}\cup B'_1$ and
$Y_{n_1}:=Y_{n_1-N}\cup B'_2$. Let $H'_n:=
\bigl(H_{n-2N}[X_{n_0-N}]\cup G_{R,1}\bigr)
+
\bigl(H_{n-2N}[Y_{n_1-N}]\cup G_{R,2}\bigr)$. Then \ref{1},~\ref{2},~\ref{3} hold. 
If $r=1$, then both $H_{n-2N}[X_{n_0-N}]$ and $G_{R,1}$ are $P_{2p_1-1}$-free, so their union is $R^-$-free. If $r\ge2$, then both graphs are $\{C_{2p_1},\ldots,C_{2p_r}\}$-free; since every member of $R^-$ contains one of these cycles, their union is $R^-$-free. 
Moreover, $H'_n[X_{n_0}]$ is nearly $(s-1)$-regular and $\max\{\Delta(H'_n[X_{n_0}]),\Delta(H'_n[Y_{n_1}])\}
\le s-1$.
When $R\cong rC_4$, $H'_n[Y_{n_1}]$ is nearly $1$-regular. 
Furthermore, $e(H'_n[B'_1\cup B'_2])=N^2+(2r-1)N/2+N/2
=g(N,N,R)$. Thus \ref{4} holds. 
For the following comparison, let $a,b\ge N+2$ be arbitrary integers with $a+b=n$.
Since $N$ is even, 
$g(a,b,R)-g(a-N,b-N,R)=N(n+r-N)$. 
When $R\not\cong rC_4$, $e(H'_n[Y_{n_1}])=1$. 
Furthermore, 
$e(H'_n[B'_1\cup B'_2])=N^2+(s-1)N/2=g(N,N,R)-1$. Thus \ref{4} holds. Since $N$ is
even, we have
$g(a,b,R)-g(a-N,b-N,R)=N(n+(s-1)/2-N)$. 
In either case, $H'_n\in\mathcal K_{n_0,n_1}(R)$ and $g(a,b,R)-g(a-N,b-N,R)$ is independent of
$(a,b)$. 
As $g(n_0-N,n_1-N,R)=f(n-2N,R)$,
Claim~\ref{claim:g-g} gives $g(n_0,n_1,R)=f(n,R)$.  Consequently, $e(H'_n)=f(n,R)$ and $H'_n\in\mathcal K_n(R)$. 
\end{proof}

\subsection{The proofs of local extremal lemmas}

\begin{proof}[\textbf{Proof of Lemma~\ref{lem:f-f}}]
By the definition of $f$ and $g$, when $R\cong 2C_6$, we have
$f(n-1,2C_6)\le
\max_{n_0+n_1=n-1}\{(n_0+3/2)(n_1+3/2)+15/4\}\le ((n+2)/2)^2+15/4$.
Moreover, we have
$f(n,2C_6)\ge g(\lceil n/2\rceil,\lfloor n/2\rfloor,2C_6)\ge
\lceil n/2\rceil\lfloor n/2\rfloor
+12\lfloor(n-1)/8\rfloor
\ge n^2/4+3n/2-49/4$.
Thus,
$f(n,2C_6)-f(n-1,2C_6)\ge n/2-17$.
The other two cases are handled similarly.
When $R\cong rC_4$, we have
$f(n-1,R)\le ((n-1+r)/2)^2-(2r-1)/4$
and
$f(n,R)\ge n^2/4+r(n-1)/2-2$.
Thus,
$f(n,R)-f(n-1,R)\ge n/2-(r^2-2r+10)/4$.
When $R\not\cong 2C_6$ and $R\not\cong rC_4$, we have
$f(n-1,R)\le
((n-1+(s-1)/2)/2)^2+1$
and
$f(n,R)
\ge n^2/4+(s-1)(n-1)/4$.
Thus,
$f(n,R)-f(n-1,R)\ge n/2-(s^2-2s+21)/16$.
\end{proof}

Before proving Lemmas~\ref{lem:special-two-C6},~\ref{lem:special-all-C4}, and~\ref{lem:bipartite}, we establish
the following facts.

\begin{fact}\label{common neighborhood}
Let $G$ be a graph with a partition
$V(G)=V_0\cup V_1$, where $|V_0|=n_0$, and write $G_i:=G[V_i]$ for $i\in\left\{0,1\right\}$. Let
$k_1:=\Delta(G_1)$, let $u\in\mathbb Z$, and let $\emptyset\ne S\subseteq V_1$. If $d_G(x)\ge n_0+u$ for  every $x\in S$, then
$|N_{G_0}(S)|
\ge n_0-|S|(k_1-u)$.
\end{fact}
\begin{proof}
For every $x\in S$,
$d_{G_0}(x)
=d_G(x)-d_{G_1}(x)
\ge n_0+u-k_1$.
Thus $x$ has at most $k_1-u$ non-neighbors in $V_0$.
The conclusion follows from the union bound.
\end{proof}

\begin{fact}
\label{fact:ex(n,tP_3)}
Let $k\ge0$ and $t\ge1$ be integers. If a graph $H$ satisfies $\Delta(H)\le k$ and contains no copy of $tP_3$, then
$e(H)\le
\left\lfloor{|V(H)|}/{2}\right\rfloor+3(t-1)k$. 
\end{fact}
\begin{proof}
Let $\mathcal P$ be a maximal family of pairwise
vertex-disjoint copies of $P_3$ in $H$, and let $T$ be the union of
their vertex sets. Since $H$ contains no copy of $tP_3$, we have
$|T|\le3(t-1)$.
By the maximality of $\mathcal P$, the graph $H-T$ is $P_3$-free,
and hence $\Delta(H-T)\le1$. Therefore,
$e(H)\le e(H-T)+\sum_{x\in T}d_H(x)
\le \left\lfloor{|V(H)|}/{2}\right\rfloor+3(t-1)k$.
\end{proof}

\begin{fact}\label{fact:pb-2p2}
Let $a,b$ be positive integers. If $H$ satisfies $\Delta(H)\le a$, contains a copy of $P_b$,
and is $P_b\cup2P_2$-free, 
then $e(H)\le ab+a+1$.
\end{fact}
\begin{proof}
Fix a copy $P$ of $P_b$ in $H$. Since $H-V(P)$ is $2P_2$-free, $\nu(H-V(P))\le1$,  so Lemma~\ref{lem:chvatal-hanson} gives $e(H-V(P))\le a+1$.
At most $ab$ edges are incident with $V(P)$,
and hence $e(H)\le ab+a+1$.
\end{proof}

\begin{proof}[\textbf{Proof of Lemma~\ref{lem:special-two-C6}}] 
Throughout the proof, we assume that $n, n_0$ are sufficiently large.
Since every graph in $\mathcal K_n(2C_6)$ has exactly
$f(n,2C_6)$ edges, one implication is immediate.
It remains to prove that $e(G)\ge f(n,2C_6)$ implies
$G\in\mathcal K_n(2C_6)$. Note that
$f(n,2C_6)=\max\{ab+6\lfloor a/4\rfloor+\binom{a\bmod4}{2}
+6\lfloor b/4\rfloor+\binom{b\bmod4}{2}:a,b\in\mathbb{Z}^+,\ a+b=n\}$.  Clearly, there exists a constant $M$ such that 
$f(n,2C_6)\ge
g(\lfloor n/2\rfloor,\lceil n/2\rceil,2C_6)
\ge n^2/4+3n/2-M$.
If $\alpha+\beta<3$, 
then for any fixed constant $M'$, 
$n_0n_1+\alpha n_0+ \beta n_1+M'\le ((n+\alpha+\beta)/2)^2+M'-\alpha\beta<f(n,2C_6)$. 
Define $V_1':=\{x\in V_1:d_G(x)\ge n_0\}$. Set $n_1':=|V_1'|$, $G_1':=G_1[V_1']$, and
$k_1':=\Delta(G_1')$. We have $k'_1\le k_1\le5$. 
If $V_1'=\emptyset$, then
$e(G)\le\lfloor5n_0/2\rfloor+(n_0-1)n_1 <f(n,2C_6)$, a contradiction.
Thus $V_1'\neq\emptyset$. Choose $x'\in V_1'$ with
$d_{G_1'}(x')=k_1'$, and let $X':=N_{G_0}\bigl(N_{G_1'}[x']\bigr)$. By Fact~\ref{common neighborhood},
$|X'|\ge n_0-k_1(k_1'+1)$. 

\noindent\textbf{Case 1.} $4\le k'_1 \le 5$.

For any $u\in V(G)$, $G[N_G(u)]$ is $(4K_1+4P_2)$-free. Otherwise, $G[N_G(u)]$ would contain a copy of $2(2K_1+2P_2)$, which would yield a copy of $\widehat{2C_6}$ in $G$, a contradiction. Thus $\nu(G_0[X'])\le3$. Note that
$\Delta(G_0[X'])\le k_0\le5$.
Lemma~\ref{lem:chvatal-hanson} gives $e(G_0[X'])\le18$. Therefore, 
\[
e(G)\le 18 +5(n_0-|X'|)+n_0n'_1+ \left\lfloor \tfrac{k'_1n'_1}{2}\right\rfloor+(n_0-1)(n_1-n'_1)
\le n_0n_1 + \left\lfloor \tfrac{5n_1}{2}\right\rfloor+O(1)<f(n,2C_6),
\]
contradicting the assumption that $e(G)\ge f(n,2C_6)$.

\noindent\textbf{Case 2.} $k'_1 =3$.

For any $u\in V(G)$, $G[N_G(u)]$
is $\big(3K_1+(P_5\cup 2P_2)\big)$-free. Otherwise, $G[N_G(u)]$ would contain a copy of $(K_1+P_5)\cup(2K_1+2P_2)$, which would yield a copy of $\widehat{2C_6}$ in $G$, a contradiction. 
Suppose first that $G_1'$ contains $P_5$.
If $G_1'$ contains a copy $F$ of $P_5\cup2P_2$, then, with $Y:=N_{G_0}(F)$, 
we have $\Delta(G_0[Y])\le2$.
By Fact~\ref{common neighborhood}, $|Y|\ge n_0-9k_1$. Hence
$e(G)\le e(G_0[Y])+5(n_0-|Y|)+n_0n_1'
+\lfloor3n_1'/2\rfloor+(n_0-1)(n_1-n_1')
<f(n,2C_6)$, a contradiction.
If $G_1'$ contains $P_5$ but is
$P_5\cup2P_2$-free, then Fact~\ref{fact:pb-2p2}, with $a=5,b=5$, gives $e(G_1')\le31$. Therefore
$e(G)\le\lfloor5n_0/2\rfloor+n_0n_1'+31
+(n_0-1)(n_1-n_1')
<f(n,2C_6)$, a contradiction. 
Thus $G_1'$ is $P_5$-free.
Since $k'_1=3$, the same observation implies that $G_0[X']$ is
$P_5\cup2P_2$-free.
If $G_0[X']$ contains a copy of $P_5$, then Fact~\ref{fact:pb-2p2} with $a=5,b=5$ yields
$e(G_0[X'])\le31$.  
Then $e(G)\le 31+5(n_0 - |X'|)+n_0n'_1+\left\lfloor {3n'_1}/{2}\right\rfloor+(n_0-1)(n_1-n'_1)
< f(n,2C_6)$,
a contradiction.
Thus $G_0[X']$ contains no $P_5$, and hence $e(G_0[X'])\le\operatorname{ex}(|X'|,P_5)$. 
Lemma~\ref{lem:ex(P_5)} yields $e(G_0[X'])\le6\left\lfloor\frac{|X'|}{4}\right\rfloor
+\binom{|X'|\bmod4}{2}\le 3|X'|/2$. 
Let $m:=n_0n'_1-e_G(V_0,V'_1)$.
Since $n_0-|X'|=O(1)$, there exists a constant $C(R)$ such that 
\[
\begin{aligned}
e(G)
&\le
\frac{3|X'|}{2}+5(n_0-|X'|)
+n_0n'_1-m+
\left\lfloor\frac{3n'_1}{2}\right\rfloor+(n_0-1)(n_1-n'_1) 
\le f(n,2C_6) + C(R)-m.
\end{aligned}
\]
Since $e(G)\ge f(n,2C_6)$, we obtain $m=O(1)$. 
We next show that $G_0$ is $P_5$-free.
If $G_0$ contains a copy $Q$ of $P_5\cup2P_2$, then let $Z:=N_{G_1'}(Q)$. Recall that  for any $u\in V(G)$, $G[N_G(u)]$ 
is $\big(3K_1+(P_5\cup 2P_2)\big)$-free. Then 
 $\Delta(G_1'[Z])\le2$.
By the definition of $m$, we have $|Z|\ge n_1'-m$, and hence
$e(G)\le {3|X'|}/{2}+5(n_0 - |X'|)+n_0n'_1+{n'_1}+3(n'_1-|Z|)+(n_0-1)(n_1-n'_1)
< f(n,2C_6)$,
a contradiction. Thus $G_0$ is $P_5\cup2P_2$-free.
If $G_0$ contains a copy of $P_5$, then Fact~\ref{fact:pb-2p2} with $a=5,b=5$ gives $e(G_0)\le31$. Consequently, $e(G)\le 31+n_0n'_1+\left\lfloor {3n'_1}/{2}\right\rfloor+(n_0-1)(n_1-n'_1)< f(n,2C_6)$, a contradiction.
Thus $G_0$ is $P_5$-free. Since $G_1'$ is also $P_5$-free, we have $e(G)\le \operatorname{ex}(n_0,P_5)+n_0n'_1-m+\operatorname{ex}(n'_1,P_5)+(n_0-1)(n_1-n'_1)
\le g(n_0,n_1,2C_6)-m-(n_1-n'_1)
\le f(n,2C_6)$ by Lemma~\ref{lem:ex(P_5)}. 
Together with $e(G)\ge f(n,2C_6)$, this forces equality throughout. In particular, $n_1=n'_1$, $m=0$, $e(G_i)=\operatorname{ex}(n_i,P_5)$ for $i\in\{0,1\}$, and
$(n_0,n_1)$ satisfies
$g(n_0,n_1,2C_6)=f(n,2C_6)$.  Therefore, $G[V_0,V_1]\cong K_{n_0,n_1}$ and, by Lemma~\ref{lem:ex(P_5)},  $G_i\cong\lfloor n_i/4\rfloor K_4\cup K_{n_i\bmod4}$ for
$i\in\{0,1\}$.
Hence $G\in\mathcal K_n(2C_6)$.

\noindent\textbf{Case 3.} $k'_1 \le 2$.

Suppose first that $G_1'$ contains a copy $T$ of $4P_2$, and let
$U:=N_{G_0}(T)$. By
Fact~\ref{common neighborhood}, $|U|\ge n_0-8k_1$.  Recall that  for any $u\in V(G)$, $G[N_G(u)]$ 
is $(4K_1+4P_2)$-free. Then $\Delta(G_0[U])\le3$.
Therefore
$e(G)\le \lfloor3|U|/2\rfloor+5(n_0-|U|)+n_0n_1'
+\lfloor k_1'n_1'/2\rfloor+(n_0-1)(n_1-n_1')
<f(n,2C_6)$, a contradiction.
Now suppose that $G_1'$ is $4P_2$-free. Since $\Delta(G_1')\le2$ and $\nu(G_1')\le3$, Lemma~\ref{lem:chvatal-hanson} gives
$e(G_1')\le9$, and hence
$e(G)\le\lfloor5n_0/2\rfloor+n_0n_1'+9
+(n_0-1)(n_1-n_1')
<f(n,2C_6)$,
a contradiction.
\end{proof}

\begin{proof}[\textbf{Proof of Lemma~\ref{lem:special-all-C4}}]
Recall that $n_0\ge n_1\ge 2r^2+r$, $n_0$ is sufficiently large, and
$g(n_0,n_1,rC_4)
=n_0n_1+\left\lfloor(2r-1)n_0/2\right\rfloor
+\left\lfloor n_1/2\right\rfloor$.
If $\alpha+\beta<2r$ and $\alpha<2r-1$, then $n_0n_1+\left\lfloor\alpha n_0/2\right\rfloor
+\left\lfloor \beta n_1/2\right\rfloor+O(1)<g(n_0,n_1,rC_4)$.
Since every graph in $\mathcal K_{n_0,n_1}(rC_4)$ has exactly
$g(n_0,n_1,rC_4)$ edges, one implication is immediate.
It remains to prove that
$e(G)\ge g(n_0,n_1,rC_4)$ implies
$G\in\mathcal K_{n_0,n_1}(rC_4)$.
Define $V_1' := \left\{x \in V_1 : d_G(x) \geq n_0\right\}$. Set $n'_1:= |V'_1|$ and $G'_1:=G_1[V'_1]$. Every vertex $x \in V_1 \setminus V_1'$ satisfies $d_G(x) \le n_0-1$.  
Thus $e(G) 
\le e(G_0)+n_0n'_1+e(G'_1)+(n_0-1)(n_1-n'_1)
\le \left\lfloor (2r-1)n_0/2\right\rfloor+n_0n_1+\left\lfloor (2r-1)n'_1/2\right\rfloor-(n_1-n'_1)$. 
If $n'_1\le 2r$, then $e(G)\le \left\lfloor (2r-1)n_0/2\right\rfloor+n_0n_1 <g(n_0,n_1,rC_4)$, a contradiction. Therefore, $n'_1 \ge 2r +1$. Set $k'_1:= \Delta(G'_1)$. 
There exists a vertex $x' \in V'_1$ such that $d_{G'_1}(x')=k'_1$. Let $X':=N_{G_0}\bigl(N_{G_1'}[x']\bigr)$. By Fact~\ref{common neighborhood}, $|X'| \ge n_0-k_1(k'_1+1)$. 

\noindent\textbf{Case 1.} $k'_1\le1$.

We have
$e(G_1')\le\lfloor n_1'/2\rfloor$, while
$e(G_0)\le\lfloor(2r-1)n_0/2\rfloor$. 
Hence
$e(G)
\le n_0n_1+\left\lfloor (2r-1)n_0/2\right\rfloor+\left\lfloor n'_1/2\right\rfloor-(n_1-n'_1)
\le g(n_0,n_1,rC_4)-(n_1-n'_1)$. 
Since $e(G)\ge g(n_0,n_1,rC_4)$ and $n'_1\le n_1$, equality holds throughout. Hence $n_1'=n_1$, $G[V_0,V_1]\cong K_{n_0,n_1}$, 
$G_1$ is nearly 1-regular and $G_0$ is nearly $(2r-1)$-regular. 
Since $e(G_1)\ge 1$ and $G$ is $\widehat{rC_4}$-free,   $G_0$ is $(rC_4)^-$-free. Hence $G \in \mathcal{K}_{n_0,n_1}(rC_4)$.

\noindent\textbf{Case 2.} $2\le k'_1\le r$.

We first show that $e(G_1')\le\lfloor n_1'/2\rfloor+3(r-1)(2r-1)$. Suppose otherwise. Then
Fact~\ref{fact:ex(n,tP_3)} yields a copy of $rP_3$ in $G_1'$.
Let $Y\subseteq V_0$ be the common neighborhood of its vertices.
By Fact~\ref{common neighborhood}, $|Y|\ge n_0-3rk_1$.
Moreover, $\Delta(G_0[Y])\le r-1$. Indeed, otherwise there exists
$y\in Y$ with $d_{G_0[Y]}(y)\ge r$. Then $G[N_G(y)]$ contains a
copy of $r(K_1+P_3)$, which would yield a copy of $\widehat{rC_4}$ in $G$, a contradiction. Using the above bounds, 
\[
\begin{aligned}
e(G)&\le \left\lfloor \frac{(r-1)|Y|}{2} \right\rfloor+(2r-1)(n_0-|Y|)+ n_0n_1+ \left\lfloor \frac{k'_1n'_1}{2} \right\rfloor -(n_1-n'_1)
<g(n_0,n_1,rC_4),
\end{aligned}
\]
a contradiction.
Thus the desired bound holds.
Since $k_1'\ge2$ and $n_1'\ge2r+1$, $G_1'$ contains a copy of $P_3\cup(2r-2)P_1$. Let $Z\subseteq V_0$ be the common neighborhood
of its vertices. By Fact~\ref{common neighborhood},
$|Z|\ge n_0-(2r+1)k_1$.
We claim that $\Delta(G_0[Z])\le2r-2$. Otherwise, there would exist 
$z\in Z$ with $d_{G_0[Z]}(z)\ge 2r-1$. Then the graph $G[N_G(z)]$ contains
a copy of $(r-1)K_{2,2}\cup(K_1+P_3)$, which would yield a copy of  $\widehat{rC_4}$ in $G$, a contradiction. Therefore, 
\[e(G)\le \left\lfloor \frac{(2r-2)|Z|}{2}\right\rfloor+(2r-1)(n_0-|Z|)+n_0n_1+\left\lfloor \frac{n'_1}{2}\right\rfloor + 3(r-1)(2r-1)-(n_1-n'_1) < g(n_0,n_1,rC_4),\]
which contradicts $e(G)\ge g(n_0,n_1,rC_4)$.

\noindent\textbf{Case 3.} $r+1\le k'_1\le 2r-2$.

Write $k_1'=2r-b$, where $2\le b\le r-1$. Recall that 
$|X'|\ge n_0-k_1(k_1'+1)\ge 2r+b$. 
We claim that $G_0[X']$ is $bP_3$-free. Otherwise,  $G[N_G(x')]$ contains
a copy of $(r-b)K_{2,2}\cup b(P_3+K_1)$, which would yield a copy of $\widehat{rC_4}$ in $G$, a contradiction.
Since $\Delta(G_0[X'])\le2r-1$,
Fact~\ref{fact:ex(n,tP_3)} gives
$e(G_0[X'])\le\lfloor|X'|/2\rfloor+3(b-1)(2r-1)$. Hence, 
\[
e(G)
\le \left\lfloor \frac{|X'|}{2}\right\rfloor+(2r-1)(3(b-1)+(n_0-|X'|))+n_0n_1+\left\lfloor \frac{k'_1n'_1}{2}\right\rfloor -(n_1-n'_1)
< g(n_0,n_1,rC_4),
\]
which contradicts $e(G)\ge g(n_0,n_1,rC_4)$.

\noindent\textbf{Case 4.} $k'_1=2r-1$.

Since the case $r=1$ was handled in Case~1, we may assume that $r\ge2$.
Recall that $|X'|\ge n_0-k_1(k'_1+1) \ge 2r+1$.
We claim that $G_0[X']$ is $P_3$-free. Otherwise, the graph $G[N_G(x')]$ contains a copy of
$(r-1)K_{2,2}\cup(P_3+K_1)$, which would yield a copy of $\widehat{rC_4}$ in $G$, a contradiction. Then $e(G_0[X'])\le \left\lfloor {|X'|}/{2}\right\rfloor $.
Let $m:=n_0n'_1-e_G(V_0,V'_1)$. There exists a constant $C(R)$ such that 
\[
\begin{aligned}
e(G) &\le \left\lfloor \frac{|X'|}{2}\right\rfloor+(2r-1)(n_0-|X'|)+n_0n_1-m+\left\lfloor \frac{(2r-1)n'_1}{2}\right\rfloor -(n_1-n'_1)\\
&\le g(n_0,n_1,rC_4)-m +C(R).
\end{aligned}
\]
Note that  $e(G)\ge g(n_0,n_1,rC_4)$. We obtain $m=O(1)$. 
If $\Delta(G_0)\le 1$, then $e(G_0)\le \left\lfloor n_0/2\right\rfloor$, and hence $e(G)
\le \left\lfloor{n_0}/{2}\right\rfloor+ n_0n_1+ \left\lfloor{(2r-1)n'_1}/{2}\right\rfloor-(n_1-n'_1)
\le  g(n_1,n_0,rC_4)-(n_1-n'_1)$.
Since $n_0\ge n_1$, $e(G)\ge 
g(n_0,n_1,rC_4)\ge g(n_1,n_0,rC_4)$.
Note that 
$r\ge 2$ and  $n_1\ge n'_1$. 
We obtain equality throughout. In particular, $n_0=n_1=n'_1$, $e(G_0)=\left\lfloor{n_0}/{2}\right\rfloor$, $e(G_1)=\left\lfloor{(2r-1)n_1}/{2}\right\rfloor$, and $G[V_0,V_1]\cong K_{n_0,n_1}$. As in Case~1, it follows that $G \in \mathcal{K}_{n_1,n_0}(rC_4)=\mathcal{K}_{n_0,n_1}(rC_4)$. 
If $\Delta(G_0)\ge2$, then $G_0$ contains a copy of $P_3\cup(2r-2)P_1$. Let $W\subseteq V_1'$ be the common neighborhood of its vertices. By the definition of $m$, we have $|W|\ge n_1'-m$.
Moreover, $\Delta(G_1'[W])\le2r-2$. 
Otherwise, there exists
$w\in W$ with $d_{G'_1[W]}(w)\ge 2r-1$. Then the graph $G[N_G(w)]$ contains a copy of $(r-1)K_{2,2}\cup(K_1+P_3)$, which would yield a copy of $\widehat{rC_4}$ in $G$, a contradiction.
Therefore, 
\[
\begin{aligned}
e(G) & \le \left\lfloor \frac{|X'|}{2}\right\rfloor
+\left\lfloor \frac{(2r-2)|W|}{2}\right\rfloor + (2r-1)(n_0-|X'|+n'_1-|W|)+n_0n_1
< g(n_0,n_1,rC_4),
\end{aligned}
\]
which contradicts $e(G)\ge g(n_0,n_1,rC_4)$.
\end{proof}

\begin{proof}[\textbf{Proof of Lemma~\ref{lem:bipartite}}]
The proof adapts Yuan's~\cite{Yuan2021} strategy for odd wheels, supplemented by additional integer and embedding lemmas needed to handle a prescribed collection of even cycles.

Throughout the proof, we assume that $n_0$ is sufficiently large.
Recall that $R\not\cong2C_6$ and $R\not\cong rC_4$, and that 
$g(n_0,n_1,R)
=n_0n_1+\left\lfloor (s-1)n_0/2\right\rfloor+1$.  
Since every graph in $\mathcal K_{n_0,n_1}(R)$ has exactly $g(n_0,n_1,R)$ edges,
one implication is immediate.
It remains to prove that
$e(G)\ge g(n_0,n_1,R)$ implies
$G\in\mathcal K_{n_0,n_1}(R)$. 
Since $G$ is $\widehat R$-free, $G[N_G(u)]$ is $R$-free for every $u\in V(G)$.
Write $\mathbf p:=(p_1,\ldots,p_r)$.
We have $p_r\ge3$, and hence $s\ge3$. We will repeatedly use the following observation: if $\alpha+\beta<s-1$ and $\beta\ge 0$, then $n_0n_1+\left\lfloor\alpha n_0/2\right\rfloor
+\left\lfloor \beta n_1/2\right\rfloor+O(1) < g(n_0,n_1,R)$. 
We first handle the case where $d_G(x)\ge n_0+1$ for all $x\in V_1$,
and then proceed by induction on $n_1$.

\begin{claim}\label{d n_0+1}
If $d_G(x)\ge n_0+1$ for every $x\in V_1$, then the lemma holds. 
\end{claim}
\begin{proof}
For every $x\in V_1$, we have
$1\le d_{G_1}(x)\le k_1$ and
$d_{G_0}(x)\ge n_0+1-k_1$.
Hence $G_1$ has no isolated vertices. 
Choose $x_0\in V_1$ with $d_{G_1}(x_0)=k_1$. 
Let $X:=N_{G_0}(N_{G_1}[x_0])$. By Fact~\ref{common neighborhood}, $|X|\ge n_0-(k_1+1)(k_1-1)$.
Fix $\varepsilon\in(0,1)$.
For $k_1\le s-2$, 
we assert that
$e(G_0[X])>
\left\lfloor{(s-k_1-2+\varepsilon)|X|}/{2}\right\rfloor$. 
Otherwise, we would have
$e(G)
\le e(G_0[X])+(n_0-|X|)(s-1)+n_0n_1+\left\lfloor{k_1n_1}/{2}\right\rfloor
< g(n_0,n_1,R)$,
a contradiction.
Thus, whenever $k_1\le s-2$, for any fixed $q$, 
Lemma~\ref{lem:many-paths,stars}, applied with
$k=s$ and $k'=s-k_1-2$, yields copies of $q P_{s-k_1}$ and $q S_{s-k_1}$ in $G_0[X]$.
We distinguish cases according to $k_1$.

\noindent\textbf{Case 1.} $k_1=1$. 

Since $G_1$ has no isolated vertices, $G_1$ is a perfect matching; moreover, $d_{G_0}(x)=n_0$ for every $x\in V_1$. 
As $n_1\ge 8s$,
the graph $G_1$ contains a copy of $2P_2\cup(s-3)P_1$.
Suppose that $k_0=s-1$, and choose $x\in V_0$ with
$d_{G_0}(x)=s-1$. Since $p_r\ge3$, the graph $G[N_G(x)]$
contains a copy of
$\bigcup_{i=1}^{r-1}K_{p_i,p_i}
\cup\bigl((p_r-1)K_1+(2P_2\cup(p_r-3)P_1)\bigr)$. Lemma~\ref{lem:path-cycle}\ref{lem:path-cycle-p2} yields a copy of $R$ in $G[N_G(x)]$, a contradiction. Hence $k_0\le s-2$.
It follows that
$e(G)\le\lfloor(s-2)n_0/2\rfloor+n_0n_1+\lfloor n_1/2\rfloor
<g(n_0,n_1,R)$,
again a contradiction.

\noindent\textbf{Case 2.} $2\le k_1\le s-3$, $(k_1,\mathbf{p}) \neq (4,(2,3,3))$ and $(k_1,\mathbf{p}) \neq (5,(3,3,3))$.

Since $\mathbf{p}\neq (3,3)$, Lemma~\ref{lem:int1} implies that either
$\sum_{i=1}^r\lceil2p_i/(s+1-k_1)\rceil\le k_1$ or
$\sum_{i=1}^r\lceil2p_i/(k_1+2)\rceil\le s-k_1$.
Suppose first that
$\sum_{i=1}^r\lceil2p_i/(s+1-k_1)\rceil\le k_1$, and set
$b_i:=\lceil2p_i/(s+1-k_1)\rceil$ for each $i\in[r]$. Then
$b_i\ge1$ and $\sum_{i=1}^r b_i\le k_1$.
Choose a copy of $\bigcup_{i=1}^r b_i\cdot P_{s-k_1}$ in $G_0[X]$ and a copy of $\bigcup_{i=1}^r b_i\cdot K_1$ in $N_{G_1}(x_0)$. 
Since all edges between $X$ and $N_{G_1}(x_0)$ are present, 
$G[N_G(x_0)]$ contains a copy of $\bigcup_{i=1}^r(b_iK_1+b_iP_{s-k_1})$.
Note that $b_i(s+1-k_1)\ge2p_i$ for each $i\in[r]$. Then  $b_i(s-k_1)\ge 2p_i-b_i$. 
Lemma~\ref{lem:path-cycle} yields a copy of $R$ in $G[N_G(x_0)]$, a contradiction.
We may therefore assume that $\sum_{i=1}^r\lceil2p_i/(s+1-k_1)\rceil>k_1$ and $\sum_{i=1}^r\lceil2p_i/(k_1+2)\rceil\le s-k_1$. Since
$(k_1,\mathbf p)\notin\{(4,(2,3,3)),(5,(3,3,3))\}$,
Lemma~\ref{lem:int2} gives
$\sum_{i=1}^r\lfloor p_i/2\rfloor\ge k_1$. Consequently,
$ s/2\le \sum_{i=1}^r\lceil p_i/2\rceil
=s-\sum_{i=1}^r\lfloor p_i/2\rfloor\le s-k_1$.

We assert that $e(G_1) \ge {3n_0}/{4}$.
Otherwise, since $e(G)\ge g(n_0,n_1,R)$, we have 
$e(G_0)\ge e(G)-e(G_1)-n_0n_1
> \left\lfloor {(s-5/2)n_0}/{2}\right\rfloor$. 
Since $k_1\ge 2$, we can fix a copy $P$ of $P_3$ in $G_1$. If $G_1$ is $(P_3\cup 2P_2)$-free, Fact~\ref{fact:pb-2p2} with $a=k_1,b=3$ yields $e(G_1)\le4k_1+1$. Since $\delta(G_1)\ge1$, $e(G_1)\ge \left\lfloor {n_1}/{2}\right\rfloor\ge4s>4k_1+1$, a contradiction. Therefore, $G_1$ contains a copy of $P_3\cup2P_2$, which extends to $P_3\cup2P_2\cup(s-5)P_1$.
Let $Y\subseteq V_0$ be the common neighborhood of its vertices. By Fact~\ref{common neighborhood}, we have $|Y|\ge n_0-(s+2)(k_1-1)$.
Since $e(G_0)
> \left\lfloor {(s-5/2)n_0}/{2}\right\rfloor$, 
$e(G_0[Y])\ge e(G_0)-(s-1)(n_0-|Y|)>\left\lfloor {(s-3)|Y|}/{2}\right\rfloor$. 
Hence there exists $y\in Y$ with $d_{G_0[Y]}(y)\ge s-2$. The graph $G[N_G(y)]$ contains a copy of $(s-2)K_1+(P_3\cup2P_2\cup(s-5)P_1)$. 
In fact, for any $u\in V(G)$, $G[N_G(u)]$ is $(s-2)K_1+(P_3\cup2P_2\cup(s-5)P_1)$-free. 
Otherwise, when $r=1$, we have $p_1=s\ge 5$, and hence $G[N_G(u)]$ contains a copy of $(p_1-2)K_1+(P_3\cup 2P_2\cup (p_1-5)P_1)$;
when $r\ge 2$, we have $p_r\ge 3$ and $p_{r-1}\ge 2$, and hence $G[N_G(u)]$ contains a copy of $\bigcup_{i=1}^{r-2} K_{p_i,p_i}\cup\big((p_{r-1}-1)K_1+(P_3\cup (p_{r-1}-2)P_1)\big)\cup \big((p_r-1)K_1+ (2P_2\cup (p_{r}-3)P_1)\big)$. 
Lemma~\ref{lem:path-cycle} with $k=p_1-2$ when $r=1$, and Lemma~\ref{lem:path-cycle}\ref{lem:path-cycle-p3} with $c=p_{r-1}-1$, together with Lemma~\ref{lem:path-cycle}\ref{lem:path-cycle-p2} when $r\ge 2$, yield a copy of $R$ in
$G[N_G(u)]$, a contradiction. Hence $e(G_1) \ge {3n_0}/{4}$, as asserted.

We now use this bound to estimate $e(G_0)$.
Recall that $s/2\le \sum_{i=1}^r\lceil p_i/2\rceil\le s-k_1\le\sum_{i=1}^rp_i$.
Hence we may choose integers $c_i'$ such that $\lceil p_i/2\rceil\le c_i'\le p_i$ for every $i\in[r]$ and
$\sum_{i=1}^r c_i'=s-k_1$.
Since $e(G_1)\ge 3n_0/4$, we have $e(G_1)\ge3n_1/4$. Moreover, $n_1\ge 2e(G_1)/k_1\ge 3n_0/(2k_1)>n_0/s$, so $n_1$ is sufficiently large. Thus Lemma~\ref{lem:many-paths,stars}, applied with $k=s$, $k'=1$ and $\varepsilon=1/2$, yields a copy of $k_1 P_3$ in $G_1$. Since $s-k_1\ge s/2$, we have $s-2k_1\ge0$.  As
$n_1\ge8s$, we can  extend the  $k_1 P_3$  to a copy of
$k_1P_3\cup(s-2k_1)P_1$ in $G_1$. Fix such a copy, and let
$Z\subseteq V_0$ be the common neighborhood of its $s+k_1$
vertices. By Fact~\ref{common neighborhood},
$|Z|\ge n_0-(s+k_1)(k_1-1)$.
We claim that $\Delta(G_0[Z])\le s-k_1-1$. Otherwise, there exists
$z\in Z$ with at least $s-k_1=\sum_{i=1}^r c_i'$ neighbors in $Z$. 
Moreover, $\sum_{i=1}^r(p_i-c_i')=k_1$ and $\sum_{i=1}^r(2c_i'-p_i)=s-2k_1$, while $2c_i'-p_i\ge0$ for every $i\in[r]$.
Hence $G[N_G(z)]$ contains a copy of $\bigcup_{i=1}^r\left(c_i'K_1+\bigl((p_i-c_i')P_3\cup(2c_i'-p_i)P_1\bigr)\right)$. Applying Lemma~\ref{lem:path-cycle}\ref{lem:path-cycle-p3} with $c=c'_i$ for each $i\in[r]$, we obtain a copy of $R$ in $G[N_G(z)]$, a contradiction.
Thus $\Delta(G_0[Z])\le s-k_1-1$. Therefore
\[
e(G_0)
\le
\left\lfloor \frac{(s-k_1-1)|Z|}{2}\right\rfloor
+(s-1)(n_0-|Z|)
\le
\left\lfloor \frac{(s-k_1-1)n_0}{2}\right\rfloor
+(s-1)(s+k_1)(k_1-1).
\]
By Lemma~\ref{lem:local-embedding-from-C1}, this implies that $G$ contains a copy of $\widehat{R}$, again a contradiction.

\noindent\textbf{Case 3.} 
$(k_1,\mathbf{p}) = (4,(2,3,3))$ or $(k_1,\mathbf{p}) = (5,(3,3,3))$.

In this case, 
$R\cong C_4\cup2C_6$ or $R\cong3C_6$, and $s=k_1+4$. 
We first show that
$e(G_1)<13n_0/8$.
Suppose that
$e(G_1)\ge13n_1/8$ and $n_1\ge n_0/2$. 
Then $n_1$ is sufficiently large.
Lemma~\ref{lem:many-paths,stars}, applied with $k=s$, $k'=3$ and
$\varepsilon=1/4$, yields a copy of $3P_5$ in $G_1$.
Let $Y\subseteq V_0$ be the common neighborhood of its vertices. By Fact~\ref{common neighborhood},
$|Y|\ge n_0-15(k_1-1)$.
For any $u\in V(G)$, $G[N_G(u)]$ is $(3K_1+3P_5)$-free.
Otherwise, $G[N_G(u)]$ contains a copy of either
$3(K_1+P_5)$ or
$(K_1+P_3)\cup2(K_1+P_5)$, which would yield a copy of $\widehat{R}$ in $G$, a contradiction.
So $\Delta(G_0[Y])\le2$ and hence $e(G) \le |Y|+15(s-1)(k_1-1)+n_0n_1+ \left\lfloor{k_1n_1}/{2} \right\rfloor< g(n_0,n_1,R)$,
a contradiction. Hence either $e(G_1) < {13n_1}/{8}$ or $n_1 < n_0/2$. In the latter case, since $k_1 \leq 5$, $e(G_1) \leq {k_1n_1}/{2} < {5n_0}/{4}$.
Consequently, in both cases, 
$e(G_1) < {13n_0}/{8}$. 
Since $k_1 \ge3$, there exists a copy of $S_4$ in $G_1$. Let $W$ be the common neighborhood of its $4$ vertices in $G_0$. By Fact~\ref{common neighborhood}, $|W| \ge n_0- 4(k_1-1)$. Thus $|W|$ is also sufficiently large. 
Since for any $u\in V(G)$, $G[N_G(u)]$ is $(3K_1+3P_5)$-free,
$G_0[W]$ is $3P_5$-free. 
By Lemma~\ref{lem:many-paths,stars} with $k=s$,  $k'=3$, $\varepsilon=1/4$, we have $e(G_0[W])<{13|W|}/{8}$. Since $s-1=k_1+3\ge7$, 
\[
e(G)
\le \frac{13}{8}|W|+4(s-1)(k_1-1)+n_0n_1+\frac{13}{8}n_0
< n_0n_1+\left\lfloor\frac{7n_0}{2}\right\rfloor+1
\le g(n_0,n_1,R),
\]
which contradicts $e(G)\ge g(n_0,n_1,R)$.

\noindent\textbf{Case 4.} $k_1=s-2$.

The case $s=3$ was handled in Case~1, so we may assume
that $s\ge4$.
We first derive an upper bound on $e(G_0)$, and then distinguish
the cases $s=4$ and $s\ge5$.
We claim that $G_0[X]$ is $\big(2P_3\cup(s-4)P_1\big)$-free. Suppose otherwise. If $r=1$, then $p_1=s\ge4$ and
$G[N_G(x_0)]$ contains a copy of
$\big(2P_3\cup(p_1-4)P_1\big)+(p_1-2)K_1$. If $r\ge2$, then
$G[N_G(x_0)]$ contains a copy of
$\bigcup_{i=1}^{r-2}K_{p_i,p_i}
\cup\bigl((P_3\cup(p_{r-1}-2)P_1)+(p_{r-1}-1)K_1\bigr)
\cup\bigl((P_3\cup(p_r-2)P_1)+(p_r-1)K_1\bigr)$.
Lemma~\ref{lem:path-cycle}\ref{lem:path-cycle-p3}, with $c=p_1-2$ when $r=1$, and with $c=p_{r-1}-1$, $c=p_r-1$ when $r\ge 2$, yields a copy of
$R$ in $G[N_G(x_0)]$, a
contradiction. 
Thus $G_0[X]$ is $\big(2P_3\cup(s-4)P_1\big)$-free. 
Recall that $|X|\ge n_0-(s-1)(s-3)\ge s+2$. It follows that $G_0[X]$ is $2P_3$-free. By Fact~\ref{fact:ex(n,tP_3)},
$e(G_0[X])\le\lfloor|X|/2\rfloor+3(s-1)$. Therefore,
$e(G_0)\le e(G_0[X])+(s-1)(n_0-|X|)
\le\lfloor n_0/2\rfloor+N(s)<3n_0/4$,
where $N(s):=(s-1)(s^2-4s+6)$. 

\noindent\textbf{Subcase 4.1.} $s=4$ and $k_1=2$. 

Since $R\not\cong2C_4$, we have $R\cong C_8$. 
In this case, $N(s)=18$.
Also,  since $e(G_0)\le 3n_0/4$,
$e(G_1)\ge e(G)-e(G_0)-n_0n_1\ge 3n_0/4 \ge 3n_1/4$.
Moreover, $e(G_1)\le n_1$, and hence $n_1\ge 3n_0/4$.
Thus $n_1$ is sufficiently large.
Lemma~\ref{lem:many-paths,stars}, applied with $k=s$, $k'=1$ and
$\varepsilon=1/2$, now yields a copy of $36P_3$ in $G_1$. Grouping these paths in pairs, we obtain $18$ pairwise vertex-disjoint copies
$Q_1,\ldots,Q_{18}$ of $2P_3$. Since $e(G_1)\le n_1$,
$e(G_0)\ge e(G)-n_0n_1-e(G_1)
\ge\lfloor n_0/2\rfloor+1$. Thus $\Delta(G_0)\ge2$. 
For any $u\in V(G)$, $G[N_G(u)]$ is $(2K_1+2P_3)$-free.
Otherwise, $G[N_G(u)]$ contains a copy of
$C_8$, a contradiction.
For every $i\in[18]$, there is a missing edge between $V(Q_i)$ and $V_0$; otherwise, choosing $u\in V_0$ with $d_{G_0}(u)\ge2$ gives a copy of $2K_1+Q_i\cong2K_1+2P_3$ in $G[N_G(u)]$. 
These missing edges are distinct, so
$e_G(V_0,V_1)\le n_0n_1-18$. Therefore
$e(G)\le\lfloor n_0/2\rfloor+18+n_0n_1-18+n_1
<g(n_0,n_1,C_8)$,
a contradiction.

\noindent\textbf{Subcase 4.2.} $s\ge 5$.

Since $e(G_0)\le\lfloor n_0/2\rfloor+N(s)<3n_0/4$, we have 
$e(G_1)\ge e(G)-e(G_0)-n_0n_1
>\lfloor(s-5/2)n_0/2\rfloor
\ge\lfloor(s-5/2)n_1/2\rfloor$.
We claim that $n_1\le n_0/2$.
Suppose otherwise. Then $n_1$ is sufficiently large. By Lemma~\ref{lem:many-paths,stars} with $k=s$, $k'=s-3$, $\varepsilon=1/2$, $G_1$ contains $N(s)\cdot S_{s-1}$. If, for each of these copies of $S_{s-1}$, not all edges between it and $V_0$ are present, then
$e(G)\le \left\lfloor{n_0}/{2}\right\rfloor+N(s)+n_0n_1-N(s)+
\left\lfloor{(s-2)n_1}/{2}\right\rfloor
<
g(n_0,n_1,R)$,
a contradiction. Hence there exists a copy $S$ of $S_{s-1}$ such that all edges between $S$ and $V_0$ are present. Recall that for any $u\in V(G)$, $G[N_G(u)]$ is $\big((s-2)K_1+(P_3\cup 2P_2\cup(s-5)P_1)\big)$-free. Then $G_0$ is $\big(P_3\cup 2P_2\cup(s-5)P_1\big)$-free. As $n_0$ is sufficiently large, $G_0$ is $\big(P_3\cup 2P_2\big)$-free. 
Since $e(G_1) \le \left\lfloor{(s-2)n_1}/{2}\right\rfloor \le \left\lfloor{(s-2)n_0}/{2}\right\rfloor$, $e(G_0)\ge e(G)-e(G_1)-n_0n_1\ge \left\lfloor{n_0}/{2}\right\rfloor+1$. 
Therefore $G_0$ contains a copy $P$ of $P_3$. Since $G_0$ is
$P_3\cup2P_2$-free, Fact~\ref{fact:pb-2p2} with $a=s-1,b=3$ yields $e(G_0)\le4s-3$, 
contradicting $e(G_0)\ge \left\lfloor{n_0}/{2}\right\rfloor+1$.  
Thus $n_1\le n_0/2$. Since $s\ge5$, 
$e(G)
\le \left\lfloor{n_0}/{2}\right\rfloor+N(s)+n_0n_1+
\left\lfloor{(s-2)n_1}/{2}\right\rfloor
\le {sn_0}/{4}+n_0n_1+N(s)
< g(n_0,n_1,R)$,  again a contradiction.

\noindent\textbf{Case 5.} $k_1=s-1$. 

Recall that $p_r\ge 3$.
For any copy $S'$ of $S_s$ with center $u$ in $G$, if $|N_G(S')|\ge s+1$, then $e(G[N_G(S')])\le 1$.
Otherwise, $G[N_G(S')]$ contains a copy of either $2P_2$ or $P_3$.
Consequently, $G[N_G(u)]$ contains a copy of either $\bigcup_{i=1}^{r-1}K_{p_i,p_i}\cup\big((p_r-1)K_1+(2P_2\cup (p_r-3)P_1)\big)$ or $\bigcup_{i=1}^{r-1}K_{p_i,p_i}\cup\big((p_r-1)K_1+(P_3\cup (p_r-2)P_1)\big)$.
In either case, Lemma~\ref{lem:path-cycle}\ref{lem:path-cycle-p3} with $c=p_r-1$ and Lemma~\ref{lem:path-cycle}\ref{lem:path-cycle-p2} yield a copy of $R$ in $G[N_G(u)]$, a contradiction. 
Note that this argument applies to every copy $S'$ of $S_s$
in any $\widehat R$-free graph $G$. We shall use this fact again
in the induction below. 
Since $|N_{G_1}(x_0)|=s-1$ and $|X|\ge n_0-s(s-2)$,
we have $e(G_0[X])\le1$. Moreover, we obtain $e(G_0)\le 1+s(s-1)(s-2)$.
If $n_1\le n_0/2$, then
$e(G)
\le O(1) + n_0n_1+\left\lfloor{(s-1)n_0}/{4}\right\rfloor
<g(n_0,n_1,R)$,
a contradiction.
Thus $n_1>n_0/2$ and $n_1$ is sufficiently large. Then
$e(G_1)\ge e(G)-e(G_0)-n_0n_1
>
\left\lfloor{(s-1)n_0}/{2}\right\rfloor-O(1)
>
\left\lfloor{(s-3/2)n_1}/{2}\right\rfloor$.
Fix an integer $M>s(s-1)(s-2)$. By Lemma~\ref{lem:many-paths,stars} with $k=s$, $k'=s-2$ and $\varepsilon=1/2$, $G_1$ contains a copy of $M S_s$.  Note that the copies are vertex-disjoint.
If every copy of $S_s$ has a missing edge to $V_0$, 
then these missing edges are distinct.
Hence $e(G)
\le 1+s(s-1)(s-2)+n_0n_1-M+
\left\lfloor{(s-1)n_1}/{2}\right\rfloor
< g(n_0,n_1,R)$,
a contradiction. 
Therefore all edges between some copy of $S_s$ and $V_0$ are present.
For any copy $S'$ of $S_s$ in $G$, if $|N_G(S')|\ge s+1$,
then $e(G[N_G(S')])\le 1$. It follows that $e(G_0)\le 1$, and hence $e(G)\le 1+n_0n_1+ \left\lfloor{(s-1)n_1}/{2}\right\rfloor$. Comparing this with $e(G)\ge 1+n_0n_1+ \left\lfloor{(s-1)n_0}/{2}\right\rfloor$ and using $n_0\ge n_1,s\ge3$,
we obtain equality throughout. In particular, $n_0=n_1$, $e(G_0)=1$, $e(G_1)=\left\lfloor(s-1)n_1/2\right\rfloor$, and $G[V_0,V_1]\cong K_{n_0,n_1}$. Since $G$ is $\widehat{R}$-free, it follows that $G_1$ is $R^-$-free. 
Thus $G \in \mathcal{K}_{n_0,n_0}(R)=\mathcal{K}_{n_0,n_1}(R)$.
\end{proof}

We now prove the lemma by induction on $n_1$. For the base case $n_1=8s$, we show that
$e(G_1)\le1$. 
Let $H:=N_{G_0}(V_1)$,
and $m:=n_0n_1-e_G(V_0,V_1)$.
Then $e(G_1)\le \left\lfloor{8s(s-1)}/{2}\right\rfloor$ and $m=e(G_1)+e(G_0)+n_0n_1-e(G)< \left\lfloor{8s(s-1)}/{2}\right\rfloor$.
Hence $|H|\ge n_0-m\ge n_0-\left\lfloor{8s(s-1)}/{2}\right\rfloor$.
If $\Delta(G_0[H])\le s-2$, then
$e(G)
\le {(s-2)|H|}/{2}+n_0n_1+O(1)
< g(n_0,n_1,R)$,
a contradiction. So $\Delta(G_0[H])= s-1$ and $G_0[H]$ contains a copy of $S_s$. Since for any copy $S'$ of $S_s$ in $G$, $|N_G(S')|\ge s+1$ implies $e(G[N_G(S')])\le 1$, we have $e(G_1)\le 1$, and hence $e(G)\le \left\lfloor{(s-1)n_0}/{2}\right\rfloor+ n_0n_1-m+1$. Combining this with $e(G)\ge \left\lfloor{(s-1)n_0}/{2}\right\rfloor+ n_0n_1+1$, we obtain equality throughout. Hence $e(G_1)=1$, $m=0$, and $e(G_0)=\left\lfloor(s-1)n_0/2\right\rfloor$. Moreover, since $G$ is $\widehat{R}$-free, $e(G_1)=1$ and
$G[V_0,V_1]\cong K_{n_0,n_1}$, the graph $G_0$ is $R^-$-free.
Hence $G\in\mathcal K_{n_0,n_1}(R)$.

Assume the result holds for $n_1-1$. We now prove it for $n_1\ge 8s+1$. 
If every $x\in V_1$ satisfies $d_G(x)\ge n_0+1$, then Claim~\ref{d n_0+1} applies.
If there exists $x\in V_1$ such that $d_G(x)\le n_0$, then
\[
e(G-x)=e(G)-d_G(x)\ge \left\lfloor\frac{(s-1)n_0}{2}\right\rfloor +n_0n_1+1-n_0=
 n_0(n_1-1)+\left\lfloor\frac{(s-1)n_0}{2}\right\rfloor+1.
\]
Since \(n_0\ge n_1-1\ge8s\), the induction hypothesis applies to \(G-x\) with the partition \(V_0\cup(V_1\setminus\{x\})\).
We have $(G-x) \in \mathcal{K}_{n_0,n_1-1}(R)$. Hence
$e(G_1-x)=1$, and $G_0$ is nearly $(s-1)$-regular and
$R^-$-free. Moreover,
$e(G-x)
=g(n_0,n_1-1,R)$ and $d_G(x)=n_0$.
It remains to show that $N_G(x)=V_0$.
Suppose otherwise that $x$ has a neighbor $y\in V_1$.
Choose $V_1'\subseteq V_1$ of size $s+1$ containing $x,y$ and
the endpoints of the unique edge of $G_1-x$.
Then $e(G_1[V_1'])\ge2$.
Let $H'=N_{G_0}(V'_1)$.
Since $G[V_0,V_1\setminus\{x\}]\cong K_{n_0,n_1-1}$, we have $H'=N_{G_0}(x)$. 
Hence $|H'|=d_{G_0}(x)\ge n_0-d_{G_1}(x)\ge n_0-(s-1)$. 
If $\Delta(G_0[H'])\le s-2$, then 
$e(G)
\le {(s-2)|H'|}/{2} + (s-1)(n_0-|H'|)
+n_0n_1+1+(s-1)
< g(n_0,n_1,R)$,
a contradiction. 
So $\Delta(G_0[H'])= s-1$ and $G_0[H']$ contains a copy of $S_s$. Since for any copy $S'$ of $S_s$ in $G$, $|N_G(S')|\ge s+1$ implies $e(G[N_G(S')])\le 1$, we have $e(G_1[V'_1])\le 1$, again a contradiction.
Thus $N_G(x)=V_0$, and hence $G\in\mathcal K_{n_0,n_1}(R)$. This completes the induction. 
\end{proof}

\section{Concluding remarks}\label{sec:concluding}

A natural direction is to replace the apex by an independent set of fixed size and study the Tur\'an problem for $tK_1+R$. Already for $t=2$ and $R=C_4$, this gives the octahedron $K_{2,2,2}$. The octahedron theorem of Erd\H{o}s and Simonovits~\cite{erdos1971extremal} states that, for sufficiently large $n$, every extremal graph for $K_{2,2,2}$ is the join of a $C_4$-extremal graph and a $P_3$-extremal graph. Thus, even this basic case connects the extremal structure to a bipartite Tur\'an problem.

For longer even cycles, a similar connection arises from a simple construction. For $k\ge3$, let $\mathcal F_k:=\{C_{2k},K_{2,k}\}$. The join of an $\mathcal F_k$-free graph and an independent set is $(2K_1+C_{2k})$-free, giving $ \operatorname{ex}(n,2K_1+C_{2k}) \ge\max_{\substack{a+b=n\\a,b\in\mathbb Z^+}}\left\{ab+\operatorname{ex}(a,\mathcal F_k)\right\}. $ Understanding the extremal structure beyond this construction is a first step towards the following general problem.

\begin{problem}
Let $t\ge2$ be a fixed integer and let $R$ be a fixed
nonempty bipartite $2$-regular graph.
Characterize the extremal graphs for $tK_1+R$
for all sufficiently large $n$.
\end{problem}

\section*{Acknowledgements}
 Junxue Zhang was supported by the National Natural Science Foundation of China (No. 12601655).

\section*{Competing interests}
The authors declare no competing interests.

\section*{Data availability}
No data were used for the research described in this article.

\section*{Declaration of generative AI use}
Generative AI (ChatGPT) was only used to assist with  language polishing. All proof ideas are entirely due to the authors, who take full responsibility for the accuracy and content of this work.
The manuscript was prepared by the authors.

\bibliographystyle{plain}
\bibliography{ref1}
\end{document}